\documentclass{elsart}
\usepackage{amssymb}

\input{tcilatex}

\begin{document}

\begin{frontmatter}

%Title, authors and addresses


\title{Tent Spaces associated with Semigroups of Operators}


\author{Tao Mei} \thanks{The author was supported in part by a Young Investigator Award of
the N.S.F supported summer workshop in Texas A\&M university 2007. }
\address{Dept. of Math. Univ. of Illinois at Urbana-Champaign, Urbana, IL, 61801, USA}

\begin{abstract}
We study tent spaces on general measure spaces $(\Omega, \mu)$. We
assume that there exists a semigroup of positive operators on
$L^p(\Omega, \mu)$ satisfying a monotone property but do not assume
any geometric/metric structure on $\Omega$. The semigroup plays the
same role as integrals on cones and cubes in Euclidean spaces. We
then study BMO spaces on general measure spaces and get an analogue
of Fefferman's $H^1$-BMO duality theory. We also get a $H^1$-BMO
duality inequality without assuming the monotone property.

All the results are proved in a more general setting, namely for
noncommutative $L^p$ spaces.

\end{abstract}


\end{frontmatter}

%\setlength{\oddsidemargin}{0in} \setlength{\evensidemargin}{0in}
%\setlength{\textwidth}{15.5cm} \setlength{\textheight}{22cm}
%\setlength{\topmargin}{-1cm}
%\begin{document}

%\begin{center}
%\mathstrut

%{\Large {\bf Tent Spaces Associated with Semigroup of operators}}

%\bigskip

%{\large Tao Mei}\footnote{
%The author was supported in part by a Young Investigator Award of the N.S.F supported summer workshop in Texas A\&M university 2007.
%}
%\end{center}

%\begin{quotation}
%{\bf Abstract.} We study Tent spaces on general measure space
%$(\Omega, \mu)$. We assume there exists a semigroup of positive
%operators on $L^p(\Omega, \mu)$ satisfying a monotone property but do not assume any
%geometric/metric structure existing on $\Omega$. The semigroup plays
%the role of integrals on cones and cubes of Euclidean space.
%We then study BMO spaces on general measure spaces and get
%an analogue of Fefferman's $H^1$-BMO duality theory. We also get a $H^1$-BMO duality inequality without assuming
%the monotone property.

%All the results are proved in a more general setting, that is for elements of a
%semifinite von Neumann algebra.
%\end{quotation}
\setcounter{section}{-1}

\section{Introduction}

\setcounter{theorem}{0}\setcounter{equation}{0} Many classical
Harmonic analysis results have been extended to more general
settings, like non Euclidean spaces, Lie groups, arbitrary measure
spaces, von Neumann algebras. We normally miss clues for such
extensions if the classical proof relies on the geometric structure
of Euclidean spaces. For examples, various integrals on cones and
cubes are used very often, as powerful techniques, in classical
analysis. But they usually do not have satisfactory analogues in the
abstract case where metric/geometric structure is not pre-defined.
However, $L^p$-spaces and semigroups of operators can be studied on
these ``domains'',
say $\Omega$, in any case. In fact, given an unbounded operator $L$ on $%
L^2(\Omega)$ with a conditionally negative kernel,
$(e^{tL})_{t\geq0}$ always provides us with a semigroup of positive
operators. It will be interesting to get an appropriate replacement
of integrals on cones and cubes by considering ``semigroup of
operators".

Tent space is a typical classical object relying on the geometric
structure of Euclidean spaces. It was introduced by Coifman, Meyer
and Stein in the 1980's (see [CMS]) and is well adapted for the
study of many subjects in classical analysis. One of the related
subjects is Fefferman's $H^1$-BMO duality theory which has been
studied in the context of semigroups by many researchers. In
particular, Varopoulos (see [V1]) established an $H^1$-BMO duality
theory for a certain family of symmetric Markovian semigroups using
a probabilistic approach. More recently, Duong/Yan studied this
topic for operators with heat kernel bounds (see [DY]). In their
proofs, the geometric structure of Euclidean spaces is essential. A
motivation of our study on tent spaces is to prove an $H^1$-BMO
duality for more general spaces.

%Another related subject
%is the equivalence between the $L^p$ norm of a function $f$ and its
%Lusin integral $ S(f)$ (see section 1.1 for its definition). An
%alternate of $S(f)$ is the Littlewood-Paley function $G(f).$ The
%equivalence of the $L^p$ norm of a function $f$ and $G(f),S(f)$ is
%usually called the Littlewood-Paley theory. In many case, the role
%of $S(f)$ is essential and can not be replaced. For example, to
%prove the $H^1$-BMO duality, one has to consider $S(f)$ if his
%achieve passes through the Carleson measures.

In this article, we define tent spaces $\mathcal{T}_p\ (p=1,\infty)$
and study the duality-relation between them for functions on
abstract domains where geometric/metric structure is unavailable. As
a replacement for the integration on cones and cubes, we consider
semigroups of positive
operators $(T_t)_t$ in the definition of our tent spaces. We prove that $%
\mathcal{T}_\infty\subset (\mathcal{T}_1)^*$ if the underlying
semigroup $(T_t)_t $ is \textit{quasi-monotone}, i.e. for some
constant $k\geq0$, $(\frac st)^k T_t-T_s$ is positive for all $s>t$
or $(\frac ts)^k T_t-T_s$ is positive for all $s<t.$ A large class
of semigroups satisfies this property. In
particular, all subordinated Poisson semigroups are quasi-monotone with $k=1$%
. We also proved that, for a quasi-monotone semigroup $(T_t)_t$, the
inverse
relation $(\mathcal{T}_1)^*\subset\mathcal{T}_\infty$ holds if and only if $%
(T_y)_y$ satisfies an \textit{$L^\frac12$ condition} :$||T_y(fT_y(g))||_{L^{%
\frac 12}}\leq c||f||_{L^1}||g||_{L^1},$ for all $y>0,f,g\geq 0$. We
prove in the appendix that a large class of semigroups (including
classical heat semigroup) on $\Bbb{R}^n$ satisfies this $L^\frac12$
condition. We have not found, unfortunately, an efficient way to
verify it for noncommutative semigroups of operators.
%Then semigroups of operators satisfy it if their kernels are bounded by heat kernels.

Using tent spaces as tools, we study $H^1$ and BMO spaces for
general semigroups of operators and get an analogue of the classical
$H^1$-BMO duality theory assuming the quasi monotone and $L^\frac12$
conditions. Without assuming these two conditions, we can only prove
a duality inequality
(see Section 4). %We also studied $H^1$ spaces
%associated with semigroups of operators without the quasi-monotone
%assumption.
%test if we get appropriate replacements or not
%We hope our study could be a good starting point of finding right
%alternates of classical techniques using integrals on cones and
%cubes.

In recent works of Junge, Le Merdy and Xu, (see [JLX], [JX]), they
consider semigroups on noncommutative $L^p$ spaces and study in
depth the corresponding maximal ergodic theory and Hardy spaces
$H^p$ for $p>1$. By using the square functions studied in [JLX],
Junge and the author obtained certain results for noncommutative
Riesz transforms in [J2] and [JM], but a
full generalization remains open. %which can be seen as a
%complementary to the previous remarkable work of Lust-Piquard's (see
%[LP1], [LP2], [LP3]).
We expect our study on general tent spaces would be helpful in the
study of noncommutative Riesz transforms since this is the case in
the classical situation. In fact, in Stein's book [St2], various
square functions are the main tools to prove the boundedness of
Riesz transforms. On the other hand, the importance of semigroups of
completely positive operators in the study of von Neumann algebras
has been impressively demonstrated due to the recent works of Popa,
Peterson and Popa/Ozawa etc. Pisier/Xu (see [PX1]) proved a $H^1$-
BMO duality for noncommutative martingales. These works motivate us
to write down the proofs of this article in the noncommutative
setting. However, it does not require much knowledge of von Neumann
algebras to
understand this paper. %As well known, all commutative
%spaces $L^p$ are special examples of noncommutative $L^p$ spaces.
%And all positive operators on commutative spaces is automatically
%completely positive.
For people whose interests are mainly the commutative case, our
proofs can be easily followed as well by regarding a von Neumann
algebra $\mathcal{M}$ as some $L^\infty(\Omega, \mu)$ and the trace
$\tau$ as a simplified notation for the integral over $\Omega$ with
respect to the measure $\mu$.

We do not assume that our semigroups admit dilations. We do not
assume they have kernels either (except in the appendix). These two
assumptions are true automatically in the classical setting but they
are not true in the general noncommutative setting.

\medskip This article is organized as follows.

Section 1 includes a brief review of classical tent spaces, basic
assumptions about the semigroup of positive operators under
consideration, definitions of our tent spaces, and a short
introduction to (noncommutative) semigroups of positive operators.
We listed our main results on tent space in Section 1.3.

In Section 2, we prove the main duality results for our tent spaces.

In Section 3, we define $H^1$ and BMO spaces associated with
semigroups and prove the desired $H^1$-BMO duality for certain
subordinated Poisson semigroups.

In Section 4, we remove the quasi-monotone assumption on the
underlying semigroup of operators and prove a duality inequality for
associated $H^1$ and BMO spaces.

%In Appendix, we verified the $L^\frac12$ condition for semigroup of
%operators on ${\Bbb R}^n$ with fast decay kernels.

\section{Preliminaries}

\subsection{Tent spaces on $\Bbb{R}\times \Bbb{R}^{+}$}

Consider a function $F:\Bbb{R}\times \Bbb{R}^{+}\longrightarrow
\Bbb{R}$. Let $A_0(F)$ be the square function defined by
\[
A_0(F)(x)=(\iint_{\Gamma _x^0}\frac 1y|F(s,y)|^2\frac{dy}yds)^{\frac 12}.
\]
where $\Gamma _x^0$ is the cone on the upper half plane with a right vertex
angle and vertex ${(x,0)}$:
\[
\Gamma _x^0=\{(s,y):|s-x|<y\}.
\]
For $1\leq p<\infty ,$ the tent space $\mathcal{T}_p$ is defined as (see
[CMS]),
\[
\mathcal{T}_p=\{F:||F||_{\mathcal{T}_p}=||A_0(F)||_{L^p}<\infty \}.
\]
Let $C(F)$ be the square function:
\[
C(F)(x)=\sup_{I}(\iint_{I\times
(0,|I|)}|F(s,y)|^2\frac{dy}yds)^{\frac 12},
\]
here the supremum is taken over all intervals $I\subset {\Bbb R}$
containing $x$. The tent space ${\mathcal{T}_\infty }$ is defined by
\[
{\mathcal{T}_\infty }=\{F:||F||_{\mathcal{T}_\infty }=||C(F)||_{L^\infty
}<\infty \}.
\]
${\mathcal{T}_\infty }$ connects to Carleson measure immediately. Recall a
Carleson measure $d\mu$ on the upper half plane is a measure satisfying
\[
\sup_{I}\iint_{I\times (0,|I|)}d\mu\leq c|I|,
\]
for all intervals $I\subset {\Bbb R}$. We see that $F\in\mathcal{T}_\infty$ if and only if the measure $d\mu=|F|^2\frac{%
dy}yds$ is a Carleson measure and
\[
||F||_{\mathcal{T}_\infty}^2=||d\mu||.
\]

A duality relation of tent spaces is proved in [CMS]. Namely
\[
\mathcal{T}_p^{*}=\mathcal{T}_q,
\]
for $1\leq p<\infty ,\frac 1p+\frac 1q=1$. %We will re-investigate
%this duality result in the case cones and cubes are unavailable.

Tent spaces have a close connection to the Hardy spaces. In fact, if
we set
\[
F(s,y)=y\nabla G(s,y)
\]
with $G$ being the harmonic extension of a function $g$ defined on $\Bbb{R}$%
, then
\begin{eqnarray*}
||F||_{\mathcal{T}_p} &\simeq& ||g||_{H^p},\ 1\leq p<\infty , \\
\text{ \quad and\ \ \ \  \quad }||F||_{ \mathcal{T}_\infty }&\simeq&
||g||_{BMO}\ (\stackrel{\mbox{def}}{=}\sup_I
(\int_I|g-g_I|^2)^{\frac 12}).
\end{eqnarray*}

The question is how to define tent spaces for general $L^p$ spaces,
for example,

$\bullet$ $L^p$ spaces on Lie groups.

$\bullet$ $L^p$ spaces on general measure spaces $(\Omega, \sigma, \mu)$.

$\bullet $ Noncommutative $L^p$ spaces.

\subsection{Semigroups of operators}

Given a measure space $(\Omega ,\sigma ,\mu )$, we consider a symmetric
diffusion semigroup of operators defined simultaneously on $L^p(\Omega )$, $%
1\leq p\leq \infty $. That is a collection of operators $(T_y)_y$ such that $%
T_{y_1}T_{y_2}=T_{y_1+y_2},T_0=id$ and

(i) $T_y$ are contractions on $L^p(\Omega)$ for all $1\leq p\leq\infty.$

(ii) $T_y$ are symmetric, i.e. $T_y=T^*_y$ on $L^2(\Omega)$.

(iii) $T_y(1)=1.$

(iv) $T_y(f)\rightarrow f$ in $L^2$ as $y\rightarrow0+$ for $f\in L^2.$

The conditions (i), (iii) above imply that the $T_y$'s are
\textit{positive operators}, i.e. $T_y(f)\geq 0$ if $f\geq 0.$

A symmetric diffusion semigroup $(T_y)$ always admits an
infinitesimal generator $L=\lim_{y\rightarrow 0}\frac {T_y-id}{y}.$
$L$ is a unbounded operator defined on $D(L)=\{f\in L^2,
\lim_{y\rightarrow 0}\frac {T_yf-f}{y}\in L^2\}$. We will write $
T_y=e^{yL}$.

The classical heat semigroup on $\Bbb{R}^n$ is a typical example of
symmetric diffusion semigroup, that is
\[
T_y=e^{y\triangle}
\]
with $\triangle =\sum_{i=1}^n\frac{\partial ^2}{\partial x_i^2}, $ the
Laplacian operator.

$T_y$ is the convolution operator with kernel $(T_yf=K_y*f)$
\begin{eqnarray}
K_y(x)=\frac{\exp (-\frac{|x|^2}{4y})}{(4\pi y)^{\frac n2}}.
\label{heat}
\end{eqnarray}

The Classical Poisson semigroup $\Bbb{R}^n$ is another popular example,
\[
P_y=e^{-y\sqrt{-\triangle }}.
\]
$P_y$ is the convolution operator with kernel
\begin{eqnarray}
K_y(x)=c_n\frac y{(|x|^2+|y|^2)^{\frac{n+1}2}}.  \label{poisson}
\end{eqnarray}

\begin{definition}
For two positive operators $T,T^{\prime }$, we write $T\geq
T^{\prime }$ if $T-T^{\prime }$ is a positive operator.
\end{definition}

By (\ref{heat}) and (\ref{poisson}), we easily see that for the classical
heat semigroup $(T_t)_t$, $T_t\leq (\frac st)^{\frac n2} T_s$ for every $t<s$%
. And for the classical Poisson semigroup $P_t$, $P_t\leq \frac ts
P_s$ for every $t>s$. Moreover, this kind of monotone property is
satisfied by all so-called subordinated Poisson semigroups.

\begin{definition}
\label{sub} Given a symmetric diffusion semigroup $(T_y)_y$ with a
generator $L$ (i.e. $T_y=e^{yL}$), the semigroup $(P_y)_y$ defined
by
\[
P_y=e^{-y\sqrt{-L}}
\]
is again a symmetric diffusion semigroup. We call it the subordinated
Poisson semigroup of $(T_y)_y$.
\end{definition}

Note $P_y$ is chosen such that
\begin{eqnarray}
(\frac{\partial ^2}{\partial y^2}+L)P_y=0.  \label{1}
\end{eqnarray}
It is well known that (see [St2])
\begin{eqnarray}
P_y=\frac 1{2\sqrt{\pi }}\int_0^\infty ye^{-\frac{y^2}{4u}}u^{-\frac
32}T_udu.  \label{idpy}
\end{eqnarray}
We can see that
\begin{eqnarray}
\frac{P_y}y(f)\downarrow \ \mathrm{as}\ y\uparrow \text{for any positive }f.
\label{sbd}
\end{eqnarray}
since $T_u$ is positive and $e^{-\frac{y^2}{4u}}u^{-\frac 32}$ is a function
decreasing with respect to $y.$ \smallskip

\subsection{Tent spaces associated with semigroups of operators}

Let $(T_y)_{y\geq 0}$ be a semigroup of operators on $L^p(\Omega ,\sigma
,\mu )$ satisfying (i)-(iv).

\begin{definition}
\label{tent} For $f\in L^2(\Omega ,L^2(\Bbb{R}_{+},\frac{dy}y))$, with $%
f_y\in L^2(\Omega )$ for each $y>0,$ we define
\begin{eqnarray*}
||f||_{\mathcal{T}_1^{(T_y)}} &=&||(\int_0^\infty T_y|f_y|^2\frac{dy}%
y)^{\frac 12}||_{L^1} \\
||f||_{\mathcal{T}_\infty ^{(T_y)}} &=&\sup_t||T_t\int_0^t|f_y|^2\frac{dy}%
y||_{L_\infty }^{\frac 12}.
\end{eqnarray*}
\end{definition}
Let $\mathcal{T}_1^{(T_y)}$ be the corresponding space after
completion. To define $\mathcal{T}_\infty ^{(T_y)}$, we need to work
a little bit more. Let $\mathcal{T}_\infty ^0=\{f\in L^2(\Omega
,L^2),||f||_{\mathcal{T}_\infty ^{(T_y)}}<\infty \}$. For a sequence
$(f^n)_n\in \mathcal{T}_\infty ^0$, we say $(f^n)_n$ $T$-converges
if $(\Pi f^n)(t)=T_t\int_0^t|f_y^n|^2\frac{dy}y$ weak-$*$ converges
in $L^\infty (\Omega )\otimes
L^\infty (\Bbb{R}_{+},dt)$. Denote this abstract limit by $%
\lim_nf^n$. Let $\mathcal{T}_\infty ^{(T_y)}$ be the space consisting of all these $%
\lim_nf^n$'s. We view $\mathcal{T}_\infty ^0$ as a subspace of
$\mathcal{T}_\infty ^{(T_y)}$ and view $\lim_nf^n$ and $\lim_ng^n$
as the same element of $\mathcal{T}_\infty ^{(T_y)}$ if
$\Pi(f^n-g^n)$ weak-$*$ converges to $0.$ Since $L^\infty (\Omega
)\otimes L^\infty (\Bbb{R}_{+})$ is weak-$*$ closed, $||\cdot ||_{\mathcal{T}%
_\infty ^{(T_y)}}$ extends to a norm on $\mathcal{T}_\infty
^{(T_y)}$ as
\[
||\lim_nf^n||_{\mathcal{T}_\infty ^{(T_y)}}=||\lim_n\Pi
(f^n)||_{L^\infty (\Omega ,L^\infty (\Bbb{R}_{+}))}^{\frac 12}.
\]
Here and in the following, $\lim_n\Pi(f^n)$ always denotes the
weak-$*$ limit of $\Pi(f^n)$.

\begin{proposition}\label{p1.1}
$\mathcal{T}_\infty ^{(T_y)}$ is complete with respect to the norm
$||\cdot ||_{\mathcal{T}_\infty ^{(T_y)}}$.
\end{proposition}
{\bf Proof.} Suppose $(f_k)_k$ is a Cauchy sequence in
$\mathcal{T}_\infty ^{(T_y)}$ with $f_k=\lim_nf_k^n,$ $f_k^n\in
\mathcal{T}_\infty ^0$. Namely, for any $\epsilon>0$,
\begin{eqnarray}
||\lim_n\Pi(f_m^n-f_j^n)||_{L^\infty (\Omega ,L^\infty
(\Bbb{R}_{+}))}<\epsilon \label{fcauchy}\end{eqnarray} for $m,j$
large enough. Since $(\lim_n\Pi f_k^n)_k$ is uniformly bounded in
$L^\infty (\Omega ,L^\infty (\Bbb{R}_{+}))$, we get a weak-$*$
convergent subsequence $(\lim_n \Pi f_{k_j}^n)_j$. Passing to the
diagonal, we get that $(\Pi f_{k_j}^j)_j$ weak-$*$ converges.
Therefore, $(f_{k_j}^j)_j$ T-converges to $f\in \mathcal{T}_\infty
^{(T_y)}.$ By (\ref{fcauchy}), for any $\epsilon
>0$,
\[
||f-f_m||_{\mathcal{T}_\infty
^{(T_y)}}=||\lim_j\Pi(f_{k_j}^j-f_{m}^j)||_{L^\infty (\Omega
,L^\infty (\Bbb{R}_{+}))}^{\frac 12}<\epsilon
\]
for $m$ large enough. This shows that $(f_k)_k$
$||\cdot||_{\mathcal{T}_\infty ^{(T_y)}}$-norm converges to $f$.\qed

Definition \ref{tent} is adapted to the classical ones because of
the following observation.

\textbf{Observation.} For a locally integrable function $f$ on $\Bbb{R}%
\times \Bbb{R}^{+},$ it is proved in [CMS] that,
\begin{eqnarray}
||A_0(f)||_{L^1}\leq ||A_k(f)||_{L^1}\leq c_{k}||A_0(f)||_{L^1},
\label{A0Ak}
\end{eqnarray}
where
\[
A_k(f)=(\iint_{\Gamma _x^k}\frac 1y|f_y|^2\frac{dy}ydt)^{\frac 12}.
\]
with $\Gamma _x^k=\{(s,y):|s-x|<2^ky\}.$ It is not hard to check
that $c_k\leq c2^{ k}$ by the ${\cal T}_1-{\cal T}_\infty$ duality
and a change of variables. We can rewrite $A_0$ and $A_k$ as square
functions of convolutions,
\begin{eqnarray*}
A_0(f) &=&(\int_0^\infty \frac 1y\chi _{(-y,y)}(\cdot )*|f_y|^2\frac{dy}%
ydt)^{\frac 12}, \\
A_k(f) &=&(\int_0^\infty \frac 1y\chi _{(-2^ky,2^ky)}(\cdot )*|f_y|^2\frac{dy%
}ydt)^{\frac 12}.
\end{eqnarray*}

If we set
\[
A_{(\mathcal{T}_y)}(f)=(\int_0^\infty T_y(|f(\cdot ,y)|^2)\frac{dy}%
ydt)^{\frac 12}.
\]
with $(T_y)_{y\geq 0}$ being a family of convolution operators with smooth
kernels $k_y$ such that
\begin{eqnarray*}
k_y(x) &>&\frac cy\ \mathrm{for}\ x\in (-y,y) \\
\mathrm{and}\ k_y(x) &\leq&  \frac
{c|y|^{1+\epsilon}}{|x|^{2+\epsilon}} \mathrm{\ as}\ |x|\rightarrow
\infty{\mathrm \ for}\ \epsilon>0,
\end{eqnarray*}
in particular, $k_y$ can be the heat kernel $K_{y^2}$, that is
\[
k_y=c\frac{\exp (-\frac{x^2}{4y^2})}y,
\]
we have
\[
cA_0^2(f)<A_{(\mathcal{T}_y)}^2(f)<\sum_k
2^{-(2+\epsilon)k}A_k^2(f).
\]
Therefore, by (\ref{A0Ak}),
\[
||A_{(\mathcal{T}_y)}(f)||_{L^p}\simeq ||A_0||_{L^p}=||f||_{\mathcal{T}_p}.
\]

We would like to search for appropriate conditions on semigroups
which provide the ``right" replacements of integrations on cones and
cubes. We pursue them by testing the duality-relation of the
associated tent spaces.

\begin{definition}
We say semigroup $(T_y)_y$ is \textbf{quasi-decreasing} if there exists ${%
{\alpha} }>0$ such that $\frac{T_y}{y^{\alpha} }$ decreases, i.e.
\begin{eqnarray}
T_t\leq (\frac ts)^{\alpha} {T_s}, \label{qdec}\end{eqnarray} for
all $0<s\leq t$.

We say $(T_y)_y$ is \textbf{quasi-increasing} if there exist
${{\alpha} }>0$ such that $y^{\alpha} T_y$ increases, i.e.
\begin{eqnarray}
 T_s\leq
(\frac st)^{\alpha} T_t, \label{qinc}
\end{eqnarray} for all $0<t\leq
s$.

We say $(T_y)_y$ is \textbf{quasi-monotone} if it is either
quasi-decreasing or quasi-increasing. \end{definition}

By (\ref{sbd}), we get

\begin{lemma}
\label{Poissondecrease} The subordinated Poisson semigroup $(P_y)_y$
of a positive semigroup $(T_y)_y$ is quasi-decreasing with
${{\alpha} }=1$.
\end{lemma}

The classical heat semigroup on $\Bbb{R}^n$ given as (\ref{heat})
satisfies the quasi-increasing condition with ${{\alpha}} =n/2.$
Heat semigroups on a complete Riemannian manifold with positive
Ricci curvature satisfy the quasi-increasing condition because of
the Harnack inequality of Li and Yau (see, for example, [Da]
Corollary 5.3.6).

We are going to prove the following duality results for our tent
spaces:

\begin{theorem}
\label{1.1} For $(T_y)_y$ quasi-monotone, we have
\[
\mathcal{T}_\infty ^{(T_y)}\subset (\mathcal{T}_1^{(T_y)})^{*}.
\]
More precisely, every $g=(g_y)_y\in \mathcal{T}_\infty ^{(T_y)}$
defines a bounded linear functional $\ell _g$ on
$\mathcal{T}_1^{(T_y)}$ by
\begin{eqnarray}
\ell _g(f)=\int_\Omega \int_0^\infty f_yg_y^{*}\frac{dy}yd\mu .
\label{l_g}
\end{eqnarray}
for all $(f_y)_y\in \mathcal{T}_1^{T_y}\bigcap L^2(\Omega ,L^2(\Bbb{R}_{+},%
\frac{dy}y))$. Here, $g_y^*$ denotes for the complex conjugate of
$g_y$, $\int_\Omega \int_0^\infty f_yg_y^{*}\frac{dy}yd\mu $
is understood as $\lim_{n} \int_\Omega \int_0^\infty f_y(g_y^{n} )^{*}%
\frac{dy}yd\mu $ for $(g_y)_y=\lim_{n} (g_y^{n} )_y$ with $(g_y^{n}
)_y\in \mathcal{T}_\infty ^0$. And
\[
||\ell _g||\leq c_{\alpha} ||(g_y)_y||_{\mathcal{T}_\infty
^{(T_y)}}.
\]
\end{theorem}

As a consequence of Theorem \ref{1.1} and Lemma
\ref{Poissondecrease}, we get

\begin{corollary}
\[
\mathcal{T}_\infty ^{(P_y)}\subset (\mathcal{T}_1^{(P_y)})^{*},
\]
with an absolute embedding constant for any subordinated Poisson semigroup $%
(P_y)_y$.
\end{corollary}

\begin{theorem}
For quasi-monotone semigroups $(T_y)_y$, we have
\begin{eqnarray}
(\mathcal{T}_1^{(T_y)})^{*}\subset \mathcal{T}_\infty ^{(T_y)},
\label{TinT2}
\end{eqnarray}
if and only if
\begin{eqnarray}
||T_y(fT_y(g))||_{L^{\frac 12}}\leq c||f||_{L^1}||g||_{L^1},\
\label{Lhalf}
\end{eqnarray}
for all $y>0,f,g\geq 0,f,g\in L^1\cap L^2.$
By (\ref{TinT2}), we mean that any linear functional $\ell $ on $\mathcal{T}%
_1^{T_y}$ is given as (\ref{l_g}) for some $g=(g_y)_y\in
\mathcal{T}_\infty ^{T_y}$ and $||(g_y)_y||_{\mathcal{T}_\infty
^{(T_y)}}\leq c_{\alpha} ||\ell ||.$
\end{theorem}

\begin{remark}{\rm We will show in the appendix that classical heat
semigroups satisfy the $L^{\frac12}$-condition (\ref{Lhalf}). And we
can see from (\ref{heat}) that
they also satisfy the quasi-monotone condition. We then get $(\mathcal{T}%
_1^{(T_y)})^*=\mathcal{T}_\infty^{(T_y)}$ for classical heat semigroups $%
(T_y)_y$. A change of variables implies $(\mathcal{T}_1^{(T_{y^2})})^*=\mathcal{%
T}_\infty^{(T_{y^2})}$. Due to the ``Observation", $\mathcal{T}%
_1^{(T_{y^2})} $'s coincide with classical tent spaces, we then
recover the duality between classical tent spaces.}
\end{remark}
As explained in the introduction, we are going to prove Theorems 1.2
and 1.4 in the noncommutative setting. We need more preliminaries
for this purpose.

\subsection{Noncommutative $L^p$ spaces and semigroups of completely
positive operators.}

Let $\mathcal{M}$ be a von Neumann algebra equipped with a normal
semifinite faithful trace $\tau $. Let $\mathcal{S}_{+}$ be the set
of all positive $x\in \mathcal{M} $ such that $\tau
(\mathrm{supp}(x))<\infty $, where $\mathrm{supp}(x)$ denotes the
support of $x$, i.e. the least projection $e\in \mathcal{M}$
such that $ex=x$. Let $\mathcal{S}_{\mathcal{M}}$ be the linear span of $%
S_{+}$. Note that $\mathcal{S}_{\mathcal{M}}$ is an involutive strongly
dense ideal of $\mathcal{M}$. For $0<p<\infty $ define
\[
\Vert x\Vert _p=\big(\tau (|x|^p)\big)^{1/p}\,,\quad x\in \mathcal{S}_{%
\mathcal{M}},
\]
where $|x|=(x^{*}x)^{1/2}$, the modulus of $x$. One can check that
$\Vert \cdot \Vert _p$ is a norm or $p$-norm on
$\mathcal{S}_{\mathcal{M}}$ according to $p\geq 1$ or $p<1$. The
corresponding completion is the noncommutative $L^p$-space
associated with $(\mathcal{M},\tau )$ and is
denoted by $L^p(\mathcal{M})$. By convention, we set $L^\infty (\mathcal{M})=%
\mathcal{M}$ equipped with the operator norm. The elements of $L^p(\mathcal{M%
})$ can be also described as measurable operators with respect to $(\mathcal{%
M},\tau )$. We refer to [PX] for more information and for more
historical references on noncommutative $L^p$-spaces. In the sequel,
unless explicitly stated otherwise, $\mathcal{M}$ will denote a
semifinite von Neumann algebra and $\tau $ a normal semifinite
faithful trace on $\mathcal{M}$.

We say an operator $T$ on $\mathcal{M}$ is completely contractive if $%
T\otimes I_n$ is contractive on $\mathcal{M}\otimes M_n$ for each
$n$. Here, $M_n$ is the algebra of $n$ by $n$ matrices and $I_n$ is
the identity operator on $M_n $. We say an operator $T$ on
$\mathcal{M}$ is completely positive if $T\otimes I_n$ is positive
on $\mathcal{M}\otimes M_n$
for each $n$. %Let $T$ be a normal contraction on ${\cal M}$. We say $T$
%is selfadjoint if
%\[
%\tau(T(x)y^*)=\tau(xT(y)^*), \ \ \forall x,y\in {\cal M}\cap L^1({\cal M}).
%\] In this case, $T$ extends to a normal contractive operator on
%$L^p({\cal M})$ for all $1\leq p\leq\infty$.

%We say a semigroup of self adjoint normal completely contractive operators $%
%(T_y)_{y\geq 0}$ on ${\cal M}$ is a noncommutative Markov semigroup
%if (i)$ T_y(1)=1$, and (iii)$T_y(f)\rightarrow f$ in the weak *
%topology of ${\cal M}$ as $y\rightarrow 0+$ for $f\in {\cal M}.$
%These conditions
%also implies $T_y$ is completely positive and $T_y$ is self adjoint on $L^2(%
%{\cal M})$. And the conditions (ii), (iii) implies $\tau T_tx=\tau
%[(T_tx)1]=\tau [x(T_t1)]=\tau[x1]=\tau x$. That means $T_t$ is trace
%preserving for any index $t$. We refer the readers to Chapter 5 of
%[JLX] for more detailed information of noncommutative Markov
%semigroups.

In this article, we will consider the so-called noncommutative
diffusion semigroup of operators $(T_y)_{y\geq 0}$ on
$L^p(\mathcal{M})$ satisfying

(i) $(T_y)_y$ are normal completely contractive on $L^p(\mathcal{M})$ for
all $1\leq p\leq \infty $.

(ii) $(T_y)_y$ are self adjoint on $L^2(\mathcal{M})$, i.e. $\tau
(T_yf)g=\tau f(T_yg),$ for all $f,g\in L^2(\mathcal{M})$.

(iii)$T_y(1)=1$,

(iv)$T_y(f)\rightarrow f$ in $L^2(\mathcal{M})$ as $y\rightarrow 0+$ for $%
f\in L^2(\mathcal{M}).$

These conditions also imply $T_y$ is completely positive and $\tau
T_tx=\tau [(T_tx)1]=\tau [x(T_t1)]=\tau[x1]=\tau x$. Namely, $T_y$'s
are trace preserving. We refer the readers to Chapter 5 of [JLX] for
more information of noncommutative diffusion semigroups.

Given a Hilbert space $H$, denote by $B(H)$ the space of all bounded
operators on $H$. Choose a norm one element $e\in H$, let $P_e$ be
the rank one projection onto Span$\{e\}.$ For $0<p\leq\infty$, let
$$L^p(\mathcal{M},H_c)=L^p(B(H)\otimes
\mathcal{M}))(1\otimes P_e).$$ Namely, $L^p(\mathcal{M},H_c)$ is the
column subspace of $L^p(B(H)\otimes \mathcal{M}))$ consisting of all
elements with the form $x(1\otimes P_e)$ for $x\in L^p(B(H)\otimes
\mathcal{M}))$. The definition of $L^p(\mathcal{M},H_c)$ does not
depend on the choice of $e$. $L^p(\mathcal{M},H_c)$ can be
identified
as the predual of $L^q(\mathcal{M},H_c)$ with $q=\frac p{p-1}$ for $%
1\leq p<\infty$. The reader can find more information on
$L^p(\mathcal{M},H_c)$ in Chapter 2 of [JLX].

All (commutative) diffusion semigroups on measurable spaces
$(\Omega,\mu)$
defined in section 1.2 are noncommutative diffusion semigroups by setting $%
\mathcal{M}=L^\infty(\Omega,\mu)$. We extend all definitions in
Section 1.2 to the noncommutative context in the natural way.

We will need the following Kadison-Schwarz inequality for unital
completely positive contraction $T$ on $L^p(\mathcal{M})$,
\begin{eqnarray}
|T(f)|^2\leq T(|f|^2),\ \ \ \ \forall f\in L^p(\mathcal{M}).
\label{cp}
\end{eqnarray}

The following definition and lemma are due to Junge/Sherman (see
[JS] Theorem 2.5).

\begin{definition}
\label{module} Let $E$ be an $\mathcal{M}$ right module with a
$L^{\frac p2}(\mathcal{M})$-valued inner product $\langle \cdot
,\cdot \rangle $. We call $E$ a Hilbert $L^p(\mathcal{M})$ $(1\leq
p<\infty )$ module if it is complete with respect to the norm
$||\cdot ||=||\langle \cdot ,\cdot \rangle ||_{L^{\frac
p2}(\mathcal{M})}^{\frac 12}.$ We call $E$ a Hilbert $ L^\infty
(\mathcal{M})$ module if it is complete with respect to the strong
operator topology generated by the seminorms
\[
||\xi ||_x=[\tau (x\langle \xi ,\xi \rangle )]^{\frac 12},x\in L^1(\mathcal{M%
}).
\]
\end{definition}

\begin{lemma}
\label{lance} Let $E$ be a Hilbert $L^p(\mathcal{M})$-module, then
$E$ is isomorphic to a complemented subspace of
$L^p(\mathcal{M},H_c)$ for some Hilbert space $H$. Moreover, the
isomorphism does not depends on $p$.
\end{lemma}

In the case of $p=\infty $, Lemma \ref{lance}  is essentially due to
Paschke (see [Pa]). The $C^{*}$-algebra analogue is due to C. Lance
(see [La], Corollary 6.3).

Let $A$ be the subspace of $L^2(\mathcal{M},L_c^2)$ such that
$a_s\in L^2({\cal M})$ for any $(a_s)_s\in A$. Define an
operator-valued inner product on the tensor product $A\otimes {\cal
M}$ by
\[
\langle (a_t)_t\otimes b,(c_t)_t\otimes d\rangle_T=b^{*}
(\int_{0}^{\infty}T_t(a_t^{*}c_t)\frac{dt}t)d.
\]
Complete $A\otimes {\cal M}$ according to Definition \ref{module} to
get a Hilbert $L^p({\cal M})$-module and denote it by
$\overline{L^{^\infty} (\mathcal{M},L_c^2)\otimes _T\mathcal{M}}^p
(p=1,\infty)$. Note the normality of $(T_s)_s$ ensures that the
inner product extends to the whole Hilbert  $L^p({\cal M})$-module.
By Lemma \ref{lance}, we get

\begin{proposition}
\label{lance2} There exist a Hilbert space $H$ and a linear map $$u:\overline{L^\infty (\mathcal{M%
},L_c^2)\otimes _T\mathcal{M}}^p\rightarrow L^p(\mathcal{M},H _c),\
\ \ p=1,\infty,$$
 such that
\[
\langle (a_t)_t\otimes b,(c_t)_t\otimes d\rangle _T=u(a\otimes
b)^{*}u(c\otimes d),
\]
for all $a\otimes b,c\otimes d\in \overline{L^\infty (\mathcal{M}
,L_c^2)\otimes _T\mathcal{M}}^p.$ And $u(\overline{L^\infty
(\mathcal{M},L_c^2)\otimes _T\mathcal{M}}^p)$ is complemented in
$L^p(\mathcal{M},H _c).$
\end{proposition}

Consider the (scalar-valued) inner product
\[
\langle a\otimes b,c\otimes d\rangle =\tau \int_0^\infty
T_s(a_s^{*}c_s)d_sb_s^{*}\frac{ds}s,
\]
for $a\otimes b,c\otimes d\in L^\infty(\mathcal{M},L_c^2)\otimes
L^2(\mathcal{M},L_c^2).$ Let $\overline
{L^\infty(\mathcal{M},L_c^2)\otimes L^2(\mathcal{M},L_c^2)}$ be the
Hilbert space completed by this inner product. We get the following
Cauchy-Schwarz inequality,
\begin{eqnarray}
|\tau \int_0^\infty a_s^{*}b_s\frac{ds}s|=|\tau \int_0^\infty a_sb_s^{*}%
\frac{ds}s| &=&|\tau \int_0^\infty a_sT_s(S_s^{-\frac 12}S_s^{\frac
12})b_s^{*}\frac{ds}s|  \nonumber \\
&\leq&[\tau \int_0^\infty T_s(S_s^{-1})|a_s|^2]^{\frac 12}[\tau
\int_0^\infty T_s(S_s)|b_s|^2]^{\frac 12},  \label{cauchy}
\end{eqnarray}
for any $(S_s)_s\geq 0$, invertible such that $(S_s^{-\frac 12}\otimes
a_s),(S_s^{\frac 12}\otimes b_s)$ are in the Hilbert space.

In this article, we will always assume our semigroup of operators
satisfy conditions (i)-(iv) listed in this section. $c_{{\alpha}} $
will be a constant depending on ${{\alpha}} $ which can be different
from line to line.
%These conditions
%can be certainly weaken. But that is not the priority of this paper.

\section{Proofs of Theorems 1.2, 1.4.}

\setcounter{theorem}{0}\setcounter{equation}{0}The noncommutative version of
Theorem 1.2 is

\begin{theorem}
\label{non1.1} If $(T_y)_{y\geq 0}$ is quasi-monotone, every $(B_s)_s\in {%
\mathcal{T}_\infty ^{(T_y)}}$ defines a bounded linear functional $\ell _B$
on ${\mathcal{T}_1^{(T_y)}}$ as
\begin{eqnarray}
\ell _B(A)={\tau }\int_0^\infty A_sB_s^{*}\frac{ds}s, \label{l_B}
\end{eqnarray}
for $(A_s)_s\in {\mathcal{T}_1^{(T_y)}}\bigcap L^2(\mathcal{M},L^2(\Bbb{R}%
_{+},\frac{dt}t)).$ And
\begin{eqnarray}
||\ell _B||\leq c_{\alpha} ||(B_s)_s||_{\mathcal{T}_\infty
^{(T_y)}}. \label{enon1.1}
\end{eqnarray}
Here, ${\tau }\int_0^\infty A_sB_s^{*}\frac{ds}s$ is understood as $\lim_{n} {%
\tau }\int_0^\infty A_s(B_s^{n} )^{*}\frac{ds}s$ for $
(B_s)_s=\lim_{n}(B_s^{n} )_s$ with $(B_s^{n} )_s\in
\mathcal{T}_\infty ^0$.
\end{theorem}

\textbf{Proof.} (i) We first prove the theorem for semigroups $(T_y)_y$
satisfying the quasi-decreasing property (\ref{qdec}) with some ${%
{\alpha}} >0$. We need the following truncated square functions $S_s,%
\widetilde{ S}_s$ in our proof:
\begin{eqnarray}
S_s &=&(\int_s^\infty T_y(|A_y|^2)\frac{y^{{{\alpha}}
-1}}{(y+s)^{{{\alpha}}} }
dy)^{\frac 12}  \label{stilde} \\
\widetilde{S}_s &=&(\int_s^\infty T_y(|A_y|^2)\frac{dy}y)^{\frac 12},
\label{stilde2}
\end{eqnarray}
for $(A_y)_y\in{\mathcal{T}_1 ^{(T_y)}}\cap L^2(\mathcal{M}, L^2(\Bbb{R}%
_+, \frac{dy}y)).$ The square functions $S_s,\widetilde{S}_s$ are
chosen to satisfy the following lemma.

\begin{lemma}
\label{lem1.1}
\begin{eqnarray}
\widetilde{S}_s &\leq &2^{\frac {\alpha} 2}S_s;  \label{lemma1} \\
\frac{dT_s(S_s)}{ds} &\geq &2T_{\frac s2}(\frac{dT_{\frac s2}(S_s)}{ds}),~%
\frac{dT_{\frac s2}(S_s)}{ds}\leq 0.  \label{lemma}
\end{eqnarray}
\end{lemma}

\textbf{Proof of Lemma \ref{lem1.1}}: (\ref{lemma1}) is obvious. We
prove (\ref{lemma}).
%To prove the first inequality of (\ref{lemma}),
Since $S_{s}\geq S_t$ for any $s\leq t$, we have
\begin{eqnarray*}
T_{s+\Delta s}(S_{s+\Delta s})-T_s(S_s) &=&T_{\frac s2}[T_{\frac {s+2\Delta
s}2}(S_{s+\Delta s})-T_{\frac s2}(S_s)] \\
&\geq&T_{\frac s2}[T_{\frac {s+2\Delta s}2}(S_{s+2\Delta s})-T_{\frac
s2}(S_s)].
\end{eqnarray*}
Divide by $\Delta s$ both sides and take $\Delta s\rightarrow 0$, we get the
first inequality of (\ref{lemma}). To prove the second inequality of (\ref
{lemma}), we apply the quasi-decreasing property of $T_s$ and get
\begin{eqnarray}
T_{y+\Delta s}(|A_y|^2)\leq T_y(|A_y|^2)(\frac{y+\Delta s}
y)^{{{\alpha}}}\leq T_y(|A_y|^2)(\frac{y+s+2\Delta s}
{y+s})^{{{\alpha}}}.  \label{tys}
\end{eqnarray}
for any $y\geq s.$ By (\ref{cp}) and (\ref{tys}), we get
\begin{eqnarray*}
&&T_{\frac{s+2\Delta s}2}S_{s+2\Delta s}-T_{\frac s2}S_s \\
&=&T_{\frac s2}T_{\Delta s}(\int_{s+2\Delta s}^\infty T_y(|A_y|^2)\frac{y^{{%
{\alpha}} -1}}{(y+s+2\Delta s)^{{{\alpha}}} }dy)^{\frac 12}-T_{\frac
s2}(\int_s^\infty T_y(|A_y|^2)\frac{y^{{{\alpha}} -1}}{(y+s)^{{{\alpha}}} }%
dy)^{\frac 12} \\
&\leq &T_{\frac s2}(\int_{s+2\Delta s}^\infty T_{y+\Delta s}(|A_y|^2)\frac{%
y^{{{\alpha}} -1}}{(y+s+2\Delta s)^{{{\alpha}}} }dy)^{\frac
12}-T_{\frac
s2}(\int_s^\infty T_y(|A_y|^2)\frac{y^{{{\alpha}} -1}}{(y+s)^{{{\alpha}}} }%
dy)^{\frac 12} \\
&\leq &T_{\frac s2}(\int_{s+2\Delta s}^\infty T_y(|A_y|^2)\frac{
y^{{{\alpha}}
-1}}{(y+s)^{{{\alpha}}} }dy)^{\frac 12}-T_{\frac s2}(\int_s^\infty T_y(|A_y|^2)%
\frac{y^{{{\alpha}} -1}}{(y+s)^{{{\alpha}}} }dy)^{\frac 12} \\
&\leq &0.
\end{eqnarray*}
Taking $\Delta s$ $\rightarrow 0$ proves the second inequality of
(\ref {lemma}).

Fix $(A_s)_s\in L^2(\mathcal{M},L_c^2)\bigcap \mathcal{T}_1^{(T_y)},(B_s)_s%
\in L^2 (\mathcal{M},L_c^2)\cap \mathcal{T}_\infty^{(T_y)}$. By
approximation, we can assume $ \widetilde{S}_s$ is invertible.  By
Lemma \ref{lem1.1} and Cauchy-Schwarz inequality (\ref{cauchy}),
\begin{eqnarray*}
|{\tau }\int_0^\infty A_sB_s^{*}\frac{ds}s| &\leq &({\tau }\int_0^\infty
T_s(|A_s|^2)\widetilde{S}_s^{-1})\frac{ds}s)^{\frac 12}({\tau }\int_0^\infty
T_s(|B_s|^2)\widetilde{S_s}\frac{ds}s)^{\frac 12} \\
&\stackrel{def}{=}&I^{\frac 12}II^{\frac 12}
\end{eqnarray*}
whenever $I,II$ are finite. Here $\widetilde{S_s}$ is defined as in (\ref
{stilde}).

For $I,$ we have
\begin{eqnarray*}
I={\tau}\int_0^\infty T_s(|A_s|^2)\widetilde{S}_s^{-1}\frac{ds}s
 ={\tau}\int_0^\infty
-\frac{d\widetilde{S}_s^2}{ds}\widetilde{S}_s^{-1}ds
=2{\tau}\int_0^\infty -\frac{d\widetilde{S}_s}{ds}ds
=2||(A_s)_s||_{\mathcal{T}_1^{(T_y)}}.
\end{eqnarray*}
For $II,$ by (\ref{lemma1}), (\ref{lemma}), we get
\begin{eqnarray*}
II \leq 2^{\frac{{\alpha}}2}{\tau}\int_0^\infty
T_s(|B_s|^2)S_s\frac{ds}s
&=&2^{\frac{{\alpha}}2}{\tau}\int_0^\infty |B_s|^2T_s(S_s)\frac{ds}s \\
&=&2^{\frac{{\alpha}}2}{\tau}\int_0^\infty |B_s|^2(-\int_s^\infty \frac{%
dT_t(S_t)}{dt}dt)\frac{ds}s \\
&=&2^{\frac{{\alpha}}2}{\tau}\int_0^\infty (\int_0^t|B_s|^2\frac{ds}s)\frac{%
dT_t(S_t)}{d(-t)}dt \\
&\leq &2\cdot 2^{\frac{{\alpha}}2}{\tau}\int_0^\infty (\int_0^t|B_s|^2\frac{ds}%
s)T_{\frac t2}(\frac{ dT_{\frac t2}(S_t)}{d(-t)})dt \\
&=&2\cdot 2^{\frac{{\alpha}}2}{\tau}\int_0^\infty T_{\frac t2}(\int_0^t|B_s|^2%
\frac{ds}s)\frac{dT_{\frac t2}(S_t)}{d(-t)}dt.
\end{eqnarray*}
Combining the estimates of I and II, we get
\begin{eqnarray*}
|{\tau}\int_0^\infty A_sB_s^{*}\frac{ds}s|\leq 4\cdot2^{\frac{{{\alpha}}}%
4}||(A_s)_s||_{\mathcal{T}_1^{(T_y)}}^{\frac 12}{\tau}\int_0^\infty T_{\frac
t2}(\int_0^t|B_s|^2\frac{ds}s)\frac{ dT_{\frac t2}(S_t)}{d(-t)}dt.
\end{eqnarray*}
Change variables and use the quasi-decreasing property of $(T_y)_y$,
we get,
\begin{eqnarray}
&&|{\tau}\int_0^\infty A_sB_s^{*}\frac{ds}s|^2  \nonumber \\
&=&|{\tau}\int_0^\infty A_{\frac s2}B_{\frac s2}^{*}\frac{ds}s|^2  \nonumber
\\
&\leq& 4\cdot2^{\frac{{{\alpha}}}2}||(A_{\frac s2})_s||_{\mathcal{T}_1^{(T_y)}}%
{\tau}\int_0^\infty T_{t}(\int_0^{\frac t2}|B_{\frac s2}|^2\frac{ds}s)\frac{%
dT_{\frac t2}(S_t)}{d(-t)}dt  \nonumber \\
&\leq& 4\cdot2^{{\alpha}}||(A_s)_s||_{\mathcal{T}_1^{(T_y)}}{\tau}%
\int_0^\infty T_{t}(\int_0^t|B_{ s}|^2\frac{ds}s)\frac{dT_{\frac
t2}(S_t)}{d(-t)}dt  \label{imply} \\
&\leq&4\cdot2^{{\alpha}}||(A_s)_s||_{\mathcal{T}_1^{(T_y)}}\sup_t||T_{t}(%
\int_0^{t}|B_s|^2\frac{ds}s)||_\infty {\tau} \int_0^\infty \frac{dT_{\frac
t2}(S_t)}{d(-t)}dt  \nonumber \\
&=&4\cdot2^{{\alpha}}||(A_s)_s||_{\mathcal{T}_1^{(T_y)}}\sup_t||T_{t}(%
\int_0^t|B_s|^2\frac{ds}s)||_\infty ||S_0||_1  \nonumber
\\
&\leq&4\cdot2^{\frac{3{\alpha}}2}||(B_s)_s||_{\mathcal{T}^{(T_y)}_\infty}^2
||(A_s)_s||_{\mathcal{T}_1^{T_y}}^2.  \label{s/2}
\end{eqnarray}
In the inequality above, we used the same notation $S_t$ for
truncated square functions of $(A_{\frac s2})_s$. Taking square root
on both sides, we proved (\ref{enon1.1}) for $(A_s)_s, (B_s)_s\in
L^2(\Omega, L^2(\Bbb{R}_+, \frac {dy}y))$ and quasi-decreasing
semigroups $(T_y)_y$. Inequality (\ref{imply}) implies that
$\lim_{n}{\tau} \int_0^\infty A_s(B_s^{n})^{*}\frac{ds}s$ exists
whenever $(B_s^{n})_s$ T-converges since $(\frac{dT_{\frac
t2}(S_t)}{d(-t)})_t\in L^1( \mathcal{M}, L^1(\Bbb{R}_+,\frac
{dt}t))$. And
\begin{eqnarray}
&&|\lim_{n}{\tau}\int_0^\infty A_s(B_s^{n})^{*}\frac{ds}s|^2
\nonumber \\
&\leq& 4\cdot2^{{\alpha}}||(A_s)_s||_{\mathcal{T}_1^{(T_y)}}{\tau}%
\int_0^\infty \lim_{n}T_{t}(\int_0^t|B_{ s}^{n}|^2\frac{ds}s)\frac{%
dT_{\frac t2}(S_t)}{d(-t)}dt  \nonumber \\
&\leq&4\cdot2^{\frac{3{\alpha}}2}||\lim_{n}(B_s^{n})_s||_{\mathcal{%
T}^{(T_y)}_\infty}^2 ||(A_s)_s||^2_{\mathcal{T}_1^{T_y}}.  \label{Tw*}
\end{eqnarray}
This means T-convergence implies weak-$*$ convergence in $(\mathcal{T}%
_1^{(T_y)})^*$. We proved Theorem (\ref{non1.1}) for quasi-decreasing
semigroups.

(ii) The proof for $(T_y)_y$ quasi-increasing requires different
truncated square functions $S_s,\widetilde{S_s}$:
\begin{eqnarray}
\widetilde{S_s} &=&(\int_s^\infty T_y(|A_y|^2)\frac{dy}y)^{\frac 12}, \\
S_s &=&(\int_s^\infty T_{2y-s}(|A_y|^2)\frac{(2y-s)^{{{\alpha}}} }{y^{{{\alpha}}%
} }\frac{ dy}y)^{\frac 12}.  \label{Ss1}
\end{eqnarray}

\begin{lemma}
\begin{eqnarray}
\widetilde{S_s} &\leq &S_s  \label{incSs} \\
\frac{dT_s(S_s)}{ds} &\geq &2T_{\frac s2}\frac{dT_{\frac s2}(S_s)}{ds},\ \
\frac{dT_{\frac s2}(S_s)}{ds}\leq 0.  \label{lemmainc}
\end{eqnarray}
\end{lemma}

\textbf{Proof. } (\ref{incSs}) is obvious by the quasi-increasing condition.
By (\ref{cp}) and the quasi-increasing condition again, it is easy to see that $%
S_s, T_{\frac s2}S_s$ are decreasing with respect to $s$. Follow the idea
used in the proof of Lemma \ref{lem1.1}, we can prove the lemma without much
difficulty.

The rest of the proof of Theorem 2.1 for quasi-increasing semigroups
is similar.\qed

We now go to prove Theorem 1.4, which is relatively easier.

The noncommutative version of Theorem 1.4 is

\begin{theorem}
\label{non1.2} Suppose semigroup $(T_y)_y$ is quasi-monotone. Then
\begin{eqnarray}
(\mathcal{T}_1^{(T_y)})^{*}\subset \mathcal{T}_\infty ^{(T_y)},
\label{T1*inTinfty}
\end{eqnarray} if and only if
\begin{eqnarray}
||T_y[(T_yg)^{\frac 12}f(T_yg)^{\frac 12}]||_{L^{\frac 12}}\leq
c||f||_{L^1}||g||_{L^1}, \label{L1/2}
\end{eqnarray}
for all $y>0,f,g\in L_{+}^1(\mathcal{M})\cap L^2(\mathcal{M})$. By (\ref
{T1*inTinfty}), we mean that any linear functional $\ell $ on $\mathcal{T}%
_1^{(T_y)}$ is given as (\ref{l_B}) for some $g=(g_y)_y\in \mathcal{T}%
_\infty ^{(T_y)}$ and
\[
||(g_y)_y||_{\mathcal{T}_\infty ^{(T_y)}}\leq c_{\alpha} ||\ell ||_{(\mathcal{T%
}_1^{(T_y)})^{*}}.
\]
\end{theorem}

\textbf{Proof.} We only prove the assertion for the quasi-increasing
case. The proof for the quasi-decreasing case is similar and
slightly easier for this Theorem. We first show that (\ref{L1/2})
implies $(\mathcal{T}_1^{(T_y)})^{*}\subseteq
\mathcal{T}_\infty ^{(T_y)}$. By Proposition \ref{lance2}, we can see $\mathcal{T}%
_1^{(T_y)}$ as a closed subspace of $L^1(\mathcal{M},H_c)$ for some
Hilbert space $H$ via the isometric embedding:
\[
f\rightarrow u(f\otimes 1).
\]
Given a linear functional $\ell \in (\mathcal{T}_1^{(T_y)})^{*},$ by
the Hahn-Banach theorem, it extends to a linear functional on
$L^1(\mathcal{M},H _c)$
with the same norm. Then there exists $\varphi \in L^\infty (\mathcal{M}%
,H_c)$ such that
\[
\ell (f)={\tau }\varphi ^{*}u(f\otimes 1).
\]
Because $u(\overline{L^\infty (\mathcal{M},L_c^2)\otimes
_T\mathcal{M}})$ is complemented in $L^\infty (\mathcal{M},H_c)$
(Proposition (\ref {lance2})), there exist
$x_n=\sum_{i=1}^na_i\otimes b_i\in {L^2(\mathcal{M},L_c^2)\otimes
\mathcal{M}}$ such that
\[
\ell (f)=\lim_{n\rightarrow \infty }{\tau }u(x_n)^{*}u(f\otimes
1)=\lim_{n\rightarrow \infty }{\tau }u(\sum_{i=1}^na_i\otimes
b_i)^{*}u(f\otimes 1),
\]
 and \[
||u(x_n)||_{L^\infty (\mathcal{M},H_c)}\leq ||\varphi ||_{L^\infty (%
\mathcal{M},H_c)}=||\ell ||.
\]
By Proposition \ref{lance2},
\begin{eqnarray}
\ell (f)=\lim_{n\rightarrow \infty }{\tau
}\sum_{i=1}^nb_i^{*}\int_{0}^{\infty}T_s(a_{i,s}^{*}f_s)\frac{ds}s
=\lim_{n\rightarrow \infty }{\tau
}\sum_{i=1}^n\int_{0}^{\infty}T_s(b_i^{*})a_{i,s}^{*}f_s\frac{ds}s.
\label{lf}
\end{eqnarray}
Set
\[
\psi _s^n=\sum_{i=1}^na_{i,s}T_{s}(b_i).
\]
It is clear that $(\psi _s^n)_s\in L^2(\mathcal{M},L_c^2)$ for each $n$ and
\[
\ell (f)=\lim_n\tau \int_0^\infty (\psi _s^n)^{*}f_s\frac{ds}s.
\]
We are going to show
\begin{eqnarray}
||T_t\int_0^t|\psi _s^n|^2\frac{ds}s||_\infty \leq
c||u(\sum_{i=0}^na_i\otimes b_i)||_{L^\infty (\mathcal{M},H_c)}.
\label{gnm}
\end{eqnarray} for $c$ independent of $t,n$. Once this is done, there
exists a subsequence of $(\psi _s^n)_s$ which T-converges to an
element $\psi \in \mathcal{T}_\infty ^{(T_y)}$ and $||\psi
||_{\mathcal{T}_\infty ^{(T_y)}}\leq c||u(\sum_{i=m}^na_i\otimes
b_i)||_{L^\infty (\mathcal{M},H_c)}\leq
c||\ell ||$ because of the weak-$*$ compactness of the unit ball of $%
L^\infty (\mathcal{M})\otimes L^\infty (\Bbb{R}_{+})$. By (\ref{Tw*}), this
will imply
\[
\ell (f)=\tau \int_0^\infty \psi _s^{*}f_s\frac{ds}s.
\]
and will prove the sufficiency of (\ref{L1/2}). We now prove
(\ref{gnm}). By the quasi-increasing property of $(T_y)_y,$ we have
\begin{eqnarray*}
||T_t\int_0^t|\psi _y^n|^2\frac{dy}y||_{L_\infty }^{\frac 12} &\leq
&2^{\frac {\alpha} 2}||T_{2t}\int_0^t|\psi
_y^n|^2\frac{dy}y||_{L_\infty
}^{\frac 12} \\
&=&2^{\frac {\alpha} 2}\sup_{{\tau }f\leq 1,f\geq 0}({\tau }%
fT_{2t}\int_0^t|\psi _y^n|^2\frac{dy}y)^{\frac 12} \\
&=&2^{\frac {\alpha} 2}\sup_{{\tau }f\leq 1,f\geq 0}({\tau }%
T_{2t}(f)\int_0^t|\psi _y^n|^2\frac{dy}y)^{\frac 12} \\
&=&2^{\frac {\alpha} 2}\sup_{{\tau }f\leq 1,f\geq 0}({\tau
}\int_0^t|\psi
_y^n(T_{2t}f)^{\frac 12}|^2\frac{dy}y)^{\frac 12} \\
&=&2^{\frac {\alpha} 2}\sup_{{\tau }f\leq 1,f\geq 0}\sup_{{\tau }%
\int_0^t|g_y|^2dy/y\leq 1}\tau \int_0^t\psi _y^n(T_{2t}f)^{\frac 12}g_y^{*}%
\frac{dy}y \\
&\leq &2^{\frac {\alpha} 2}\sup_f\sup_{g_y}||(\psi _y^n)_y||_{(\mathcal{T}%
_1)^{*}}||(g_y(T_{2t}f)^{\frac 12})_{0<y<t}||_{\mathcal{T}_1}.
\end{eqnarray*}
Note in the inequality above, we can restrict the supremum to be
taken for $f,g$ very nice, so that $||(g_y(T_{2t}f)^{\frac
12})_{0<y<t}||_{\mathcal{T}_1}$ make sense. By (\ref{lf}), we have
\[
||(\psi _y^n)_y||_{(\mathcal{T}_1)^{*}}\leq
||u(\sum_{i=0}^na_i\otimes b_i)||_{L^\infty (\mathcal{M},H_c)}.
\]
Therefore,
\begin{eqnarray}
&&||T_{t}\int_0^{t}|\psi _y^n|^2\frac{dy}y||_{L_\infty }\nonumber\\
&\leq& 2^{\alpha}||u(\sum_{i=0}^na_i\otimes b_i)||_{L^\infty
(\mathcal{M} ,H_c)}^2\sup_{f,g_y}||(g_y(T_{2t}f)^{\frac
12})_{0<y<t}||_{\mathcal{T}_1}^2 . \label{supsup}
\end{eqnarray}
Apply the Kadison-Schwarz inequality, we get
\begin{eqnarray*}
\tau [\int_{0}^{t}T_y|g_y(T_{2t}f)^{\frac 12}|^2\frac{dy}y]^{\frac
12}
&=&{\tau }T_{t}[\int_{0}^{t}T_y|g_y(T_{2t}f)^{\frac 12}|^2\frac{dy}%
y]^{\frac 12} \\
&\leq &{\tau }[\int_{0}^{t}T_{y+t}|g_y(T_{2t}f)^{\frac 12}|^2\frac{dy}%
y]^{\frac 12} \\
&\leq &2^{\frac {\alpha} 2}{\tau }[T_{2t}\int_{0}^{t}|g_y(T_{2t}f)^{\frac 12}|^2\frac{dy}%
y]^{\frac 12} \\
&=&2^{\frac {\alpha} 2}{\tau }\left( T_{2t}[(T_{2t}f)^{\frac
12}\int_0^{t}|g_y|^2\frac{dy}y(T_{2t}f)^{\frac 12}]\right) ^{\frac
12}.
\end{eqnarray*}
Using (\ref{L1/2}) for $g=\int_0^{t}|g_y|^2\frac{dy}y$, we get
\begin{eqnarray}
\tau [\int_{0}^{t}T_y|g_y(T_{2t}f)^{\frac 12}|^2\frac{dy}y]^{\frac
12}\leq c2^{\frac {\alpha} 2}||f||_{L^1}^{\frac
12}||\int_0^{t}|g_y|^2\frac{dy}y||_{L^1}^{\frac 12}. \label{p11}
\end{eqnarray} Combine (\ref{p11}), (\ref{supsup}) and take the supremum
over $k$, we get
\[
||(\psi _y^n)_y||_{\mathcal{T}_\infty ^{(T_y)}}\leq
||u(\sum_{i=0}^na_i\otimes b_i)||_{L^\infty (\mathcal{M},H_c)},
\]
which is (\ref{gnm}). We then proved the sufficiency of (\ref{L1/2}).

To prove the necessity of (\ref{L1/2}), we are going to show the necessity
of the following stronger inequality
\begin{eqnarray*}
{\tau }(\int_0^tT_y[(T_tf)^{\frac 12}|g_y|^2(T_tf)^{\frac 12}]\frac{dy}%
y)^{\frac 12}\leq ||\int_0^t|g_y|^2\frac{dy}y||_1^{\frac 12}||f||_1^{\frac
12},
\end{eqnarray*}
for all $f\in L^1_+\cap L^2_+,g\in L^1(\mathcal{M}, L^2_c)\cap L^2(\mathcal{M%
}, L^2_c).$ To see it is stronger than (\ref{L1/2}), one can consider $g_y=%
\sqrt{t}g^{\frac 12}\frac{\chi _{(t-\epsilon ,t)}(y)}{\sqrt{\epsilon
}}$ and send $\epsilon \rightarrow 0$. Assume that $
(\mathcal{T}_1^{(T_y)})^{*}\subset \mathcal{T}_\infty ^{(T_y)}$. Fix
$f$ , $(g_y)_y$, we have
\begin{eqnarray*}
{\tau }(\int_0^tT_y[(T_tf)^{\frac 12}|g_y|^2(T_tf)^{\frac 12}]\frac{dy}%
y)^{\frac 12} &=&||(g_y(T_tf)^{\frac 12})_{0<y<t}||_{\mathcal{T}_1^{(T_y)}}
\\
&\leq &\sup_{||(h_y)_y||_{\mathcal{T}_\infty \leq 1}}{\tau }%
\int_0^tg_y(T_tf)^{\frac 12}(h_y)^{*}\frac{dy}y \\
&\leq &[{\tau }\int_0^t|g_y|^2\frac{dy}y]^{\frac 12}\sup_{||(h_y)_y||_{%
\mathcal{T}_\infty \leq 1}}[{\tau }\int_0^t(T_tf)|h_y|^2\frac{dy}%
y]^{\frac 12} \\
&=&||\int_0^t|g_y|^2\frac{dy}y||_1^{\frac 12}\sup_{||(h_y)_y||_{%
\mathcal{T}_\infty \leq 1}}[{\tau }\int_0^tfT_t|h_y|^2\frac{dy}%
y]^{\frac 12} \\
&\leq &||\int_0^t|g_y|^2\frac{dy}y||_1^{\frac 12}||f||_1^{\frac
12}\sup_{||(h_y)_y||_{\mathcal{T}_\infty \leq
1}}||\int_0^tT_t|h_y|^2\frac{dy}y||_\infty ^{\frac 12} \\
&\leq &||\int_0^t|g_y|^2\frac{dy}y||_1^{\frac 12}||f||_1^{\frac 12}.
\end{eqnarray*}
The proof of the theorem is complete.\qed

\begin{remark}{\rm
\label{weakmonotone} From the proof,  we see that the quasi-monotone
assumption in Theorem \ref{non1.2} can be replaced by a ``weaker''
condition: $T_{2s}\leq cT_s,$ for all $s$ or $T_s\leq cT_{2s},$ for
all $s$. }\end{remark}

\begin{remark}{\rm\label{equivalent}
Applying the same technique used in the proof of Theorem
\ref{non1.2}, it is not hard to show that the noncommutative
$L^{\frac 12}$ condition (\ref {L1/2}) is equivalent to any of the
following conditions:

(i)$||T_t|h|^2||_{L^\infty({\cal M})}^{\frac 12}\leq c\sup_{\tau
(T_t|f|^2)^{\frac 12}\leq 1}|\tau fh^*|$, for any $t>0, h\in
L^2({\cal M}).$

(ii) $||T_t|\sum_{k=1}^n T_t(b_k)a_k|^2||_{\infty}\leq
c||\sum_{k,j=1}^nb_k^*T_t(a_k^*a_j)b_j||_{\infty}$ for any
$n\in{\Bbb N},(a_k)_{k=1}^n$, $(b_k)_{k=1}^n\in L^\infty({\cal M})$.
}\end{remark}

\begin{remark}{\rm
\label{change} By changing variables $y\rightarrow y^2$ and setting
$A_y^{\prime
}=A_{y^2},B_y^{\prime }=B_{y^2}$, we see that the duality between $\mathcal{T}%
_1^{(T_y)}$ and $\mathcal{T}_\infty ^{(T_y)}$ holds if and only if the
duality between $\mathcal{T}_1^{(T_{y^2})}$ and $\mathcal{T}_\infty
^{(T_{y^2})}$ holds. Let $(T_y)_y$ be the classical heat semigroup defined
as in (\ref{heat}), ``Observation'' in Section 1.3 tells us that $\mathcal{T}%
_p^{(T_{y^2})}$ coincide with the classical ones. We recover the
duality between the classical $\mathcal{T}_1$ and
$\mathcal{T}_\infty $ by Theorems 1.2, 1.4 (or Theorems 2.1, 2.3)
since the classical heat semigroup is quasi-increasing and satisfies
the $L^{\frac 12}$ condition (\ref{L1/2}) (see a proof in the
appendix). }\end{remark}

We will need the following results in Section 3.
\begin{lemma}
\label{T2s} Suppose a semigroup $(T_y)_y$ is quasi-monotone and satisfies the $%
L^{\frac 12}$ condition (\ref{L1/2}). We have
\[
||(T_{2s}A_s)_s||_{\mathcal{T}_1^{(T_s)}}^2\leq c_{\alpha} ||(A_s)_s||_{%
\mathcal{T}_1^{(T_s)}}\tau (\int_0^\infty |T_sA_s|^2\frac{ds}s)^{\frac 12}.
\]
\end{lemma}

\textbf{Proof.} The assumption of the lemma implies the duality between $%
\mathcal{T}_1^{(T_s)}$ and $\mathcal{T}_\infty^{(T_s)}$, which yields that
\begin{eqnarray*}
||(T_{2s}A_s)_s||_{\mathcal{T}_1^{(T_s)}}&\leq&c_{{{\alpha}}%
}\sup_{||(B_s)_s||_{\mathcal{T}_\infty^{(T_s)}}\leq1} \tau\int_0^\infty
T_{2s}(A_s)B_s\frac{ds}s.
\end{eqnarray*}
We now estimate $\tau\int_0^\infty T_{2s}(A_s)B_s\frac{ds}s$
following the proof of Theorem \ref{non1.1}. We will benefit because
of the extra $T_{2s}$. Let $S_s,\widetilde{S}_s$ be as in the proof
of Theorem \ref {non1.1} and set
\[
G_s=(\int_s^\infty|T_yA_y|^2\frac{dy}y)^\frac12.
\]
Then $G_s\leq \widetilde{S}_s\leq 2^{\frac {{{\alpha}}} 2}S_s$. By
the Cauchy-Schwarz inequality, we have
\begin{eqnarray*}
|\tau\int_0^\infty T_{2s}(A_s) B_s\frac{ds}s| &=&|\tau\int_0^\infty T_s(A_s)
T_s(B_s)\frac{ds}s| \\
&\leq&(\tau\int_0^\infty |T_s(A_s)|^2G_s^{-1}\frac{ds}s)^{\frac12}
(\tau\int_0^\infty |T_s(B_s)|^2G_s\frac{ds}s)^{\frac12} \\
&\leq&2^{\frac {{{\alpha}}} 4}(\tau\int_0^\infty |T_s(A_s)|^2G_s^{-1}\frac{ds}%
s)^{\frac12} (\tau\int_0^\infty |T_s(B_s)|^2S_s\frac{ds}s)^{\frac12} \\
&\leq&2^{\frac {{{\alpha}}} 4}(\tau\int_0^\infty |T_s(A_s)|^2G_s^{-1}\frac{ds}%
s)^{\frac12} (\tau\int_0^\infty T_s|B_s|^2S_s\frac{ds}s)^{\frac12} \\
&\stackrel{def}{=}&2^{\frac {{{\alpha}}} 4}I^{\frac 12}II^{\frac 12}
\end{eqnarray*}
We get exactly the same ``II'' as in the proof of Theorem \ref{non1.1}. Then
\[
II\leq c||(B_s)_s||_{\mathcal{T}_\infty^{(T_s)}}||(A_s)||_{\mathcal{T}%
_1^{(T_s)}}.
\]
And
\begin{eqnarray*}
I=\tau\int_0^\infty |T_s(A_s)|^2G_s^{-1}\frac{ds}s
&=&\tau\int_0^\infty -\frac{\partial G_s^2}{\partial s}G_s^{-1}ds \\
&=&-\tau\int_0^\infty \frac{\partial G_s}{\partial s}G_sG_s^{-1}+G_s\frac{%
\partial G_s}{\partial s}G_s^{-1}ds \\
&=&-2\tau\int_0^\infty \frac{\partial G_s}{\partial s}ds =2\tau G_0
\end{eqnarray*}
Therefore,
\begin{eqnarray*}
\tau\int_0^\infty T_{2s}(A_s) B_s\frac{ds}s\leq c_{{{\alpha}}%
}||(B_s)_s||^{\frac12}_{\mathcal{T}_\infty^{(T_s)}}||(A_s)_s||_{\mathcal{T}%
_1^{(T_s)}}^{\frac12}(\tau G_0)^{\frac12}.
\end{eqnarray*}
Taking the supremum over $(B_s)_s$ we get
\[
||(T_{2s}A_s)_s||_{\mathcal{T}_1^{(T_s)}}^2\leq c_{{{\alpha}}} ||(A_s)_s||_{%
\mathcal{T}_1^{(T_s)}}\tau(\int_0^\infty |T_sA_s|^2\frac{ds}s)^{\frac12}.%
\qed
\]

\begin{proposition}
\label{alpha} Assume $(T_y)_y$ is quasi monotone and satisfies the $L^{\frac
12}$ condition (\ref{L1/2}). Then for any family $(A_s)_{s\geq 0},$
\[
||(A_s)_s||_{\mathcal{T}_1^{(T_s)}}\stackrel{c_{{{{\alpha}} }}}{\thickapprox }%
||(A_s)_s||_{\mathcal{T}_1^{(T_{2s})}}.
\]
\end{proposition}

\textbf{Proof.} For $T_{s}$ quasi-increasing, we have for any $(A_s)_s$, $%
(B_s)_s$,
\begin{eqnarray}
||(A_{s})_s||_{\mathcal{T}_1^{(T_s)}}\leq 2^{\frac {{{\alpha}}%
}2}||(A_{s})_s||_{\mathcal{T}_1^{(T_{2s})}}, \ \ \ \ ||(B_{s})_s||_{\mathcal{%
T}_\infty^{(T_s)}}\leq 2^{\frac {{{\alpha}}}2}||(B_{s})_s||_{\mathcal{T}%
_\infty^{(T_{2s})}}.  \label{Bs}
\end{eqnarray}
Note the assumption of the lemma implies the duality between $\mathcal{T}%
_1^{(T_{2s})}$ and $\mathcal{T}_\infty^{(T_{2s})}$. This duality and (\ref{Bs}%
) yield that
\[
||(A_s)_s||_{\mathcal{T}_1^{(T_{s})}}\stackrel{ c_{{{\alpha}} }}{\thickapprox }%
||(A_{s})_s||_{\mathcal{T}_1^{(T_{2s})}}.
\]

The proof for quasi-deceasing $(T_s)_s$ is similar. % by the following
%2 facts:
%\begin{eqnarray*}
%||(A_{s})_s||_{{\cal T}_1^{(T_s)}} &\leq&2^{\frac
%{{{\alpha}}}2}||(A_{2s})_s||_{{\cal T}_1^{(T_s)}}.
%\end{eqnarray*}
%and
%\begin{eqnarray*}
%||(B_{s})_s||_{{\cal T}_\infty^{(T_s)}}&=&
%||\sup_t \int_0^{2t} T_{2t}|B_{s}|^2\frac{ds}s||_\infty^{\frac 12}\\
%&\leq& 2^{\frac{{{\alpha}}}2}||\sup_t \int_0^{2t} T_{t}|B_s|^2\frac{ds}s||_\infty^{\frac 12}\\
%&\leq& 2^{\frac{{{\alpha}}}2}||\sup_t \int_0^{t} T_{t}|B_{2u}|^2\frac{du}u||_\infty^{\frac 12}\\
%&=&2^{\frac{{{\alpha}}}2}||(B_{2s})_s||_{{\cal T}_\infty^{(T_s)}}.\qed
%\end{eqnarray*}
\qed

\section{$H^1-$BMO duality for Subordinated Poisson semigroups}

\setcounter{theorem}{0}\setcounter{equation}{0}

Consider the subordinated Poisson Semigroup $(P_y)_y$ of a symmetric
diffusion semigroup $(T_y)_y$. We are going to study BMO spaces
associated with $(P_y)_y.$ We first define a seminorm for
$\varphi\in L^2({\cal M})$ as
\[
||\varphi ||_{BMO_c(P)}=\sup_{y>0}||P_y(|\varphi -P_y\varphi
|^2)||_\infty ^{\frac 12}.
\]

For a sequence $(\varphi_{n} )_{n}\in L^2({\cal M}),$ with
$||\varphi_n||_{BMO_c(P)}<\infty$, let $\Phi_{n}$ be the operator
valued function $\Phi_{n}(y)=P_y(|\varphi_{n} -P_y\varphi_{n} |^2)$.
We say $(\varphi_{n} )_{n}$ P-converges if $(\Phi_{n})_{n}$ weak-$*$
converges in $L^\infty ({\cal M} )\otimes L^\infty
(\Bbb{R}_{+},dy)$. Denote this abstract limit of $(\varphi_{n} )_{n}
$ by $\lim_{n} \varphi_{n} $. Add $\lim_{n} \varphi_{n} $'s to
$\{\varphi\in L^2({\cal M}), {||\varphi||_{BMO_c(P)}<\infty}\}$ and
denote the new vector space by BMO$_c(P)$. Since the weak-$*$ limit
of $(\Phi_{n})_{n}$ exists in $L^\infty ({\cal M})\otimes L^\infty
(\Bbb{R}_{+},dy)$, $||\cdot ||_{BMO_c(P)}$ extends to a seminorm on
BMO$_c(P)$ as
\[
||\lim_{n} \varphi_{n} ||_{BMO_c(P)}=||\lim_{n}P_y(|\varphi_{n}
-P_y\varphi_{n} |^2)||_{L^\infty ({\cal M})\otimes L^\infty
(\Bbb{R}_{+})}^\frac12.
\]

Similar to Proposition \ref{p1.1}, BMO$_c(P)$ is complete with
respect to the seminorm ${||\cdot||}_{BMO_c(P)}$ because the unit
ball of $L^\infty ({\cal M})\otimes L^\infty (\Bbb{R}_{+})$ is
weak-$*$ compact. We view BMO$_c(P)$ as the resulting Banach space
after quotienting out {$\{||\varphi||_{{\rm BMO}_c(P)}=0\}$.}

In the classical case (i.e for functions $\varphi$ on ${\Bbb R}$),
it is well known that $||\varphi||_{BMO}\approx \sup_{z\in {\Bbb
R}\times {\Bbb R}^+}P_z|\varphi-P_z\varphi|$ with $P_z$ the Poisson
integral at the point $z$ (see [Ga] P217, [Pe] P79). Our definition
of BMO is an analogue of this characterization. The difference is
that $P_z\varphi$ is a number while $P_y\varphi$ is a function. And
$P_z|\varphi-P_z\varphi|\neq P_y|\varphi-P_y\varphi|(x)$ for
$z=(x,y)$ in general.

In [JM], we proved that BMO$_c(P)$ (combining with the row space)
serves as an end point of $L^p({\cal M})$ for interpolation. The
goal of this section is to find an $H^1$ space as the predual of
$BMO_c(P)$. The main tool will be the duality result of our tent
spaces in Section 2. So we need first prove a relation between
BMO$_c(P)$ and ${\cal T}_\infty^{(P_y)}$.

Let $\Gamma $ be the gradient form associated with the generator
$L,$ i.e.
\begin{eqnarray}
2\Gamma (x,y)=L(x^{*}y)-L(x^{*})y-x^{*}L(y).  \label{carre}
\end{eqnarray}

Let $\widetilde{\Gamma }$ be the gradient form associated with the
new generator $\widetilde{L}=L+\frac{\partial ^2}{
\partial s^2}$ defined on a dense subset of $L^2({\cal M}\otimes L^\infty({\Bbb R}_+))$.
Namely, $2\widetilde\Gamma (x,y)=\widetilde L(x^{*}y)-\widetilde
L(x^{*})y-x^{*}\widetilde L(y).$ By the definition, we get
\begin{eqnarray} \widetilde{\Gamma }(x,y)=\Gamma (x,y)+\frac
\partial {\partial s}x^{*}\frac
\partial {\partial s}y. \label{Gamma>partial}
\end{eqnarray}

\begin{proposition} For any $x\in L^2({\cal M})$,
\begin{eqnarray}
\Gamma (x,x) &\geq &0  \label{carrepositive} \\
\widetilde{\Gamma }(P_sx,P_sy) &=&\widetilde{L}((P_sx)^{*}P_sy).
\label{carrL}
\end{eqnarray}
\end{proposition}

{\bf Proof.} (\ref{carrepositive}) can be proved by considering the
derivative of $ e^{sL}(|e^{(t-s)L}x|^2)$ with respect to $s$ and
letting $t,s\rightarrow 0.$ In fact, $ \frac {\partial
e^{sL}(|e^{(t-s)L}x|^2)}{\partial s}=
e^{sL}\Gamma(e^{(t-s)L}x,e^{(t-s)L}x)$. (\ref{carrL}) can be seen by
the fact $\widetilde{L}(P_sx)=0$ for all $x$.\qed

\medskip

%For a pair of element $f\in L^1({\cal M}),\varphi \in {\cal M},$ we
%can set
%\begin{eqnarray*}
%A_s &=&\frac{\partial P_s}{\partial s}f,B_s=\frac{\partial P_s}{\partial s}%
%\varphi \\
%\mbox{or }A_sB_s^{*} &=&\Gamma (P_sf,P_s\varphi ),\mbox{ or
%}A_sB_s^{*}=\widetilde{\Gamma }(P_sf,P_s\varphi ).
%\end{eqnarray*}
%We first give an analogue of the classical relation between Tent spaces $%
%{\cal T}_\infty $ and the BMO spaces.

\begin{theorem} For any $\varphi\in L^2({\cal M})$, we have
\label{tbmo}
\begin{eqnarray}
||(s\frac{\partial P_s}{\partial s}(\varphi-P_{s}\varphi))_s
||_{{\cal T}_\infty ^{(P_s)}}\leq c||\varphi
||_{BMO_c(P)}.\label{carlbmo}
\end{eqnarray}
Moreover,  if $(\varphi_{n})_{n}\in L^2({\cal M})\cap BMO_c(P)$
P-converges then $(s\frac{\partial P_s}{\partial
s}(\varphi_{n}-P_{s}\varphi_{n}))_{n}$ T-converges in ${\cal
T}_\infty ^{(P_s)}$ and
\begin{eqnarray}
||\lim_{n}(s\frac{\partial P_s}{\partial
s}(\varphi_{n}-P_{s}\varphi_{n}))_s ||_{{\cal T}_\infty
^{(P_s)}}\leq c||\lim_{n}\varphi_{n} ||_{BMO_c(P)}.\label{carlbmo1}
\end{eqnarray}
\end{theorem}

{\bf Convention}. Because of Theorem \ref{tbmo} we understand
$(s\frac{\partial P_s}{\partial s}(\varphi-P_s\varphi))_s$ as an
element in ${{\cal T}_\infty ^{(P_s)}}$ via the corresponding
T-limit for any $\varphi\in BMO_c(P)$.

To prove Theorem \ref{tbmo} we need the following Lemma.

\begin{lemma}\label{subharmonic}
For any $y\geq 0, \varphi\in L^2({\cal M})$, we have
\begin{eqnarray}
\int_0^\infty P_{s+y}\widetilde{\Gamma }(P_{s}\varphi ,P_{s}\varphi )\frac{%
sy}{s+y}ds\leq P_y(|\varphi |^2). \label{Gammatild}
\end{eqnarray}
\end{lemma}

{\bf Proof}. Fix a scaler $y$ and a positive element $z\in
L^\infty({\cal M}),$ (\ref{carrL}) implies
\begin{eqnarray*}
&&{\tau (}z\int_0^\infty P_{s+y}\widetilde{\Gamma }(P_{s}\varphi
,P_{s}\varphi )\frac{sy}{s+y}ds) \\
&=&{\tau }\int P_{s+y}(z)\frac{sy}{s+y}\widetilde{L}(|P_{s}\varphi |^2)ds \\
&=&{\tau }\int L(P_{s+y}(z)\frac{sy}{s+y})|P_{s}\varphi |^2ds+{\tau
}\int (P_{s+y}(z)\frac{sy}{s+y})\frac{\partial ^2}{\partial
s^2}|P_{s}\varphi |^2ds.
\end{eqnarray*}
We use integration by parts to the second term and get
\begin{eqnarray*}
&&{\tau }\int_0^\infty (P_{s+y}(z)\frac{sy}{s+y})\frac{\partial
^2}{\partial
s^2}|P_{s}\varphi |^2ds\\
&=&0-{\tau }\int \frac
\partial {\partial s}(P_{s+y}(z)\frac{sy}{s+y})\frac \partial
{\partial s}|P_{s}\varphi |^2ds \\
&=&{\tau }\int \frac{\partial ^2}{\partial
s^2}(P_{s+y}(z)\frac{sy}{s+y})|P_{s}\varphi |^2ds+{\tau
}P_{0+y}(z)|P_0\varphi |^2\\
&=&{\tau }\int (\frac{\partial ^2}{\partial
s^2}P_{s+y}(z))\frac{sy}{s+y}|P_{s}\varphi |^2ds+{%
\tau }\int [2\frac \partial {\partial s}P_{s+y}(z)\frac \partial {\partial s}%
\frac{sy}{s+y}\\&&+P_{s+y}(z)\frac{\partial ^2}{\partial
s^2}\frac{sy}{s+y}]|P_{s}\varphi |^2ds+{\tau }P_y(z)|\varphi |^2.
\end{eqnarray*}
In the process of integration by parts above, we used the fact
$\frac
\partial {\partial s}P_s\varphi=0$ as $s=\infty$, which can be seen from the inequality (\ref{dps=0}) below.
 Thus, by the definition of
$\widetilde{L}$, we have
\begin{eqnarray*}
&&{\tau (}z\int_0^\infty P_{s+y}\widetilde{\Gamma }(P_{s}\varphi
,P_{s}\varphi )\frac{sy}{s+y}ds)\\
&=&{\tau }\int \widetilde{L}(P_{s+y}(z))\frac{sy}{s+y}|P_{s}\varphi
|^2ds+{
\tau }\int [2\frac \partial {\partial s}P_{s+y}(z)\frac \partial {\partial s}%
\frac{sy}{s+y} \\
&&+P_{s+y}(z)\frac{\partial ^2}{\partial
s^2}\frac{sy}{s+y}]|P_{s}\varphi
|^2ds+{\tau }P_y(z)|\varphi |^2 \\
&=&0+{\tau }\int [2\frac \partial {\partial s}P_{s+y}(z)\frac
\partial
{\partial s}\frac{sy}{s+y}+P_{s+y}(z)\frac{\partial ^2}{\partial s^2}\frac{sy%
}{s+y}]|P_{s}\varphi |^2ds+{\tau }P_y(z)|\varphi |^2 \\
&=&{\tau }\int \frac{2y^2}{(s+y)^2}[\frac \partial {\partial
s}P_{s+y}(z)-\frac 1{s+y}P_{s+y}(z)]|P_{s}\varphi |^2ds+{\tau }%
P_y(z)|\varphi |^2.
\end{eqnarray*}
By (\ref{sbd}), we have $\frac \partial {\partial s}(\frac{P_{s+y}(z)}{s+y}%
)\leq 0.$ That is
\begin{eqnarray}
\frac{\partial P_{s+y}(z)}{\partial s}\frac 1{s+y}-\frac
1{(s+y)^2}P_{s+y}(z)\leq 0. \label{dps=0}
\end{eqnarray}
Then
\[
{\tau }\int \frac{2y^2}{(s+y)^2}[\frac \partial {\partial
s}P_{s+y}(z)-\frac 1{s+y}P_{s+y}(z)]|P_{s}\varphi |^2ds\leq 0.
\]
Therefore
\[
\tau(z\int_0^\infty P_{s+y}\widetilde{\Gamma }(P_{s}\varphi ,P_{s}\varphi )\frac{%
sy}{s+y}ds)\leq {\tau }P_y(z)|\varphi |^2=\tau (zP_y|\varphi |^2).
\]
By the arbitrariness of $z,$ we proved the Lemma.\qed

\smallskip
{\bf Proof of Theorem \ref{tbmo}.} Given a $\varphi\in L^2({\cal
M})$, we split $\frac{\partial P_s}{\partial s}(\varphi -P_s\varphi
)$ into three parts
\begin{eqnarray*}
\frac{\partial P_s}{\partial s}(\varphi -P_s\varphi ) &=&\frac{\partial P_s}{%
\partial s}(\varphi -P_y\varphi )+\frac{\partial P_s}{\partial s}%
(P_{s+y}\varphi -P_y\varphi )+\frac{\partial P_s}{\partial
s}(P_s\varphi
-P_{s+y}\varphi ) \\
&=&A+B+C.
\end{eqnarray*}
It is easy to derive from (\ref{idpy}) and (\ref{cp}) that $|\frac{\partial P_y}{\partial y}(x)|^2\leq c%
\frac{P_{\frac y2}}{y^2}|x|^2$. Apply this property to B, we get
\begin{eqnarray}
P_y\int_0^y|B|^2sds&=&P_y\int_0^y|\frac{\partial P_y}{\partial y}%
P_s(P_s\varphi -\varphi )|^2sds \nonumber\\
&\leq &\frac cyP_y\int_0^yP_{\frac y2}P_s|P_s\varphi -\varphi |^2ds \nonumber\\
&=&\frac cyP_{\frac{3y}2}\int_0^yP_s|P_{s}\varphi -\varphi |^2ds.
\label{B}
\end{eqnarray}
For the terms A, C,  by (\ref{Gamma>partial}), we have $|\frac{\partial P_s}{\partial s}%
\varphi |^2\leq \widetilde{\Gamma }(P_s\varphi ,P_s\varphi ).$ Then,
by (\ref {Gammatild}) and (\ref{cp}), we get
\begin{eqnarray}
P_y\int_0^y|A|^2sds&\leq& 2P_y|\varphi -P_y\varphi |^2.\label{A}\\
P_y\int_0^y|C|^2sds&\leq& \int_0^yP_{y+s}|A|^2sds\leq 2P_y|\varphi
-P_y\varphi |^2. \label{C}
\end{eqnarray} Combine the estimates of A, B, C, we get, for any
$\varphi\in L^2({\cal M})$,
\[
P_y\int_0^y|\frac{\partial P_s}{\partial s}(\varphi -P_s\varphi
)|^2sds\leq c||\varphi ||_{BMO_c(P)}.
\]
On the other hand, by (\ref{B}), for any $f(y)\in
L^1(\mathcal{M}\otimes L^\infty (\Bbb{R}_{+},dy)),$ we have
\begin{eqnarray}
\tau \int_0^\infty (P_y\int_0^y|B|^2sds)f(y)dy=\tau \int_0^\infty
P_s|P_{s}\varphi -\varphi |^2(\int_s^\infty \frac
cyP_{\frac{3y}2}f(y)dy)ds. \label{fy}
\end{eqnarray}
Since $\int_s^\infty \frac cyP_{\frac{3y}2}f(y)dy\in
L^1(\mathcal{M}\otimes L^\infty (\Bbb{R}_{+},ds)),$ we conclude from
(\ref{fy}), (\ref{A}) and (\ref{C}) that $s\frac{\partial
P_s}{\partial s}(\varphi_{n} -P_s\varphi_{n} )$ T-converges in
${\cal T}_\infty^{P_s}$ if $\varphi_{n}$ P-converges in BMO$_c(P)$
and
\[||\lim_{n}(s\frac{\partial P_s}{\partial
s}(\varphi_{n}-P_{s}\varphi_{n}))_s||_{{\cal T}_\infty ^{(P_s)}}\leq
c||\lim_{n}\varphi_{n} ||_{BMO_c(P)}.\qed
\]

%\begin{remark}{\rm
%Assume the so-called ``$\Gamma^2\geq0$" condition: $\Gamma (P_sx,P_sx)\leq P_s\Gamma(x,x),$ we can replace $|\frac{\partial P_s}{%
%\partial s}\varphi |^2$ in (\ref{carlbmo}) with $\widetilde{\Gamma }%
%(P_s\varphi ,P_s\varphi ).$
%}\end{remark}

As an immediate consequence of Theorems \ref{non1.1} and \ref{tbmo},
we get
\begin{corollary} \label{suggest}
For any subordinated Poisson semigroup $(P_y)_y$, we have
\begin{eqnarray}
|{\tau }f\varphi ^{*}|\leq c||(\int_0^\infty P_y|\frac {\partial
P_yf}{\partial y}|^2ydy)^\frac12||_{L^1}||\varphi ||_{BMO_c(P)},
\end{eqnarray}
for any $\varphi\in L^2({\cal M}), f\in L^2({\cal M}).$
\end{corollary}

{\bf Proof.}  We know from (\ref{sbd}) that any subordinated Poisson
semigroup $(P_y)_y$ is quasi decreasing with ${{\alpha} }=1$.
Applying Theorems \ref{non1.1} and \ref{tbmo}, we get
\begin{eqnarray*}
|\tau \varphi^*f| &=&9|\tau\int_0^\infty [\frac{\partial P_y%
}{\partial y}(\varphi-P_y\varphi) ^{*}\frac{%
\partial P_y}{\partial y}f]ydy|\\
&\leq &c||(y\frac{\partial P_y}{\partial y}f)_y||_{%
\mathcal{T}_1^{(P_y)}}||(y\frac{\partial P_y}{\partial y}[\varphi
-P_{y}\varphi ])_y||_{\mathcal{T}_\infty ^{(P_y)}} \\
&\leq&c||(\int_0^\infty P_y|\frac{\partial P_yf}{\partial
y}|^2ydy)^{\frac 12}||_{L^1}||\varphi ||_{BMO_c(P)}.\qed
\end{eqnarray*}

Corollary \ref{suggest} suggests an $H^1$ norm of $f$:
$||(\int_0^\infty P_y|\frac {\partial P_yf}{\partial
y}|^2ydy)^\frac12||_{L^1}$. However, this norm does not fit the
classical case. In fact, if $P_y$ is the classical Poisson integral
operator on ${\Bbb R}^n$, $||(\int_0^\infty P_y|\frac {\partial
P_yf}{\partial y}|^2ydy)^\frac12||_{L^p}$ is equivalent to
$||f||_{H^p({\Bbb R}^n)}$ only when $p>\frac {n+1}2$. We have to
consider a smaller norm for general $H^1$ if we want to cover the
classical case.

Consider the tent space $\mathcal{T}_1^{(T_{y^2})}$ associated with ($%
T_{y^2})_{y\geq 0}$. Remark \ref{change} explains that the duality
result for $\mathcal{T}_1^{(T_{y})}$ applies to
$\mathcal{T}_1^{(T_{y^2})}$. Given $f\in L^2(\mathcal{M})$, it is
easy to see that $(y\frac{\partial P_y}{\partial y}f)_y\in
L^2(\mathcal{M},L_c^2)$. We say that $f$ belongs to the Hardy space
$H_c^1(P)$ if $(y\frac{\partial P_y}{\partial y}f)_y$ belongs to
$\mathcal{T}_1^{(T_{y^2})}$. Set
\[
||f||_{H_c^1(P)}=||(y\frac{\partial P_y}{\partial y}f)_y||_{\mathcal{T}%
_1^{(T_{y^2})}}.
\]
An equivalent definition is
\[
||f||_{H_c^1(P)}=||S(f)||_{L^1}
\]
with
\[
S(f)=(\int_0^\infty T_{y^2}|\frac{\partial P_y}{\partial
y}f|^2ydy)^{\frac 12}.
\]
Let $H_c^1(P)$ be the corresponding space after completion.
$H_c^1(P)$ can be viewed as a closed subspace of
$\mathcal{T}_1^{(T_{y^2})}$ via the embedding: $f\mapsto
(y\frac{\partial P_y}{\partial y}f)_y$.

We will show that
\[
BMO_c(P)\subseteq (H_c^1(P))^{*}
\]
provided ($T_y)_{y\geq 0}$ is quasi monotone. And
\[
BMO_c(P)=(H_c^1(P))^{*}
\]
if ($T_y)_{y\geq 0}$ satisfies the $L^{\frac 12}$ condition
(\ref{L1/2}) too.

\begin{theorem}\label{3.5}
Assume the underlying semigroup $(T_y)_y$ is quasi-monotone. Then
BMO$_c(P) \subseteq(H_c^1(P))^{*}.$ More precisely, every
$\varphi\in BMO_c(P)$ defines a linear functional $\ell_{\varphi}$
on $H_c^1(P)$ by $ \ell_{\varphi}(f)={\tau }f\varphi ^{*},$ for any
$f\in H_c^1(P)\cap L^2({\cal M})$.  And
\begin{eqnarray}
|\ell_{\varphi}|\leq c||\varphi ||_{BMO_c(P)}. \label{h1bmo}
\end{eqnarray}
Here ${\tau }f\varphi ^{*}$ is understood as $\lim_{n}{\tau
}f\varphi_{n} ^{*}$ for $\varphi$ being a P-limit of
$(\varphi_{n})_{n}\in L^2({\cal M}).$
\end{theorem}

{\bf Proof. } By the identity (\ref{idpy}), for ${(T_y)}_y$
quasi-increasing, we have,
\begin{eqnarray}
P_y=\frac 1{2\sqrt{\pi }}\int_0^\infty ye^{-\frac{y^2}{4u}}u^{-\frac
32}T_udu  \geq \frac 1{2\sqrt{\pi
}}\int_{y^2}^{2y^2}ye^{-\frac{y^2}{4u}}u^{-\frac 32}T_udu \geq
cT_{y^2}.  \label{Py>Ty2}
\end{eqnarray}

For ${(T_y)}_y$ quasi-decreasing,
\begin{eqnarray}
P_y=\frac 1{2\sqrt{\pi }}\int_0^\infty ye^{-\frac{y^2}{4u}}u^{-\frac
32}T_udu  \geq \frac 1{2\sqrt{\pi
}}\int_{\frac{y^2}2}^{y^2}ye^{-\frac{y^2}{4u} }u^{-\frac
32}T_udu\geq cT_{y^2}.  \label{Py>Ty22}
\end{eqnarray}

(\ref{Py>Ty2}), (\ref{Py>Ty22}) and Theorem \ref{tbmo} imply that
$(y\frac{\partial P_y}{\partial y}(\varphi-P_y\varphi))_y\in{{\cal
T}_\infty ^{T_{y^2}}}$ for $\varphi\in {BMO_c(P)}\cap L^2({\cal M})$
and
\begin{eqnarray}
||(y\frac{\partial P_y}{\partial y}(\varphi-P_y\varphi))_y ||_{{\cal
T}_\infty ^{(T_{y^2})}}\leq c||\varphi ||_{BMO_c(P)}.  \label{last}
\end{eqnarray}

Combining (\ref{last}) and Remark \ref{change} we get
\begin{eqnarray*}
|{\tau }f\varphi ^{*}| &=&9|\tau\int_0^\infty (y\frac{\partial
P_y}{\partial y}f)y
\frac{\partial P_y}{\partial y}(\varphi-P_y\varphi)^{*}\frac{dy}y| \\
&\leq &c||(y\frac{\partial P_y}{\partial y}f)_y||_{{\cal T}_1
^{(T_{y^2})}}||(y
\frac{\partial P_y}{\partial y}(\varphi-P_y\varphi))_y ||_{{\cal T}_\infty ^{(T_{y^2})}} \\
&\leq &c||S(f)||_{L^1}||\varphi
||_{BMO_c(P)}=c||f||_{H_c^1(P)}||\varphi ||_{BMO_c(P)}.
\end{eqnarray*}

By Theorem \ref{tbmo} and the end of the proof (i) of Theorem
\ref{non1.2}, we see that $\lim_{n}{\tau }f\varphi_{n} ^{*}$ is well
defined for $f\in L^2({\cal M})\cap H_c^1(P)$ and a P-convergent
sequence $(\varphi_{n})_{n}$. Moreover,
$$|\lim_{n}{\tau }f\varphi_{n} ^{*}|\leq c||f||_{H_c^1(P)}||\varphi ||_{BMO_c(P)}.$$
This proves Theorem \ref{3.5}.\qed

We now go to show the other direction of the desired duality result.
In the classical case, this direction is relatively easier. But it
is really complicated in our case due to the missing of the
geometric structure on von Neumann algebras (in particular, the
general measure spaces).

\begin{proposition}For $(T_y)_y$ quasi-monotone, $\varphi\in
L^2({\cal M})$, we have
\[
||\varphi ||_{BMO_c(P)}\thickapprox ||\sup_tT_{t^2}|\varphi
-P_t\varphi |^2||_\infty ^{\frac 12}.
\]
%Let $H_G^p$ be the Hardy space characterized by the norm
%\[
%||f||_{H_G^p}=||G(f)||_{L^p}
%\]
%with
%\[
%G(f)=(\int_0^\infty |\frac{\partial P_y}{\partial y}f|^2ydy)^{\frac 12}.
%\]
\end{proposition}

{\bf Proof.} By (\ref{Py>Ty2}) and (\ref{Py>Ty22}), we have
$$T_{y^2}(f)\leq c_{\alpha} P_y(f)$$
for any positive $f$. Then
$$||\sup_tT_{t^2}|\varphi -P_t\varphi
|^2||_\infty ^{\frac 12}\leq c_{\alpha}||\varphi ||_{BMO_c(P)}.
$$
On the other hand, by the identity (\ref{idpy})
\begin{eqnarray*}
 ||P_{t}|\varphi -P_t\varphi
|^2||_\infty^{\frac12}&=&||\frac 1{2\sqrt{\pi }}\int_0^\infty
te^{-\frac{t^2}{4u}}u^{-\frac 32}T_u|\varphi -P_t\varphi
|^2du||_\infty^{\frac12}\\
&\leq&||\frac 1{2\sqrt{\pi }}\int_0^{t^2}
te^{-\frac{t^2}{4u}}u^{-\frac
32}T_u|\varphi -P_t\varphi|^2du||_\infty^{\frac12}\\
&&+||\frac 1{2\sqrt{\pi }}\int_{t^2}^\infty
te^{-\frac{t^2}{4u}}u^{-\frac
32}T_u|\varphi -P_t\varphi|^2du||_\infty^{\frac12}\\
&\leq&(\frac 1{2\sqrt{\pi }}\int_0^{t^2}
te^{-\frac{t^2}{4u}}u^{-\frac
32}||T_u|\varphi -P_t\varphi|^2||_\infty du)^{\frac12}\\
&&+(\frac 1{2\sqrt{\pi }}\int_{t^2}^\infty
te^{-\frac{t^2}{4u}}u^{-\frac 32}||T_u|\varphi
-P_t\varphi|^2||_\infty du)^{\frac12}.
\end{eqnarray*}
Note for $u\geq t^2$,
\begin{eqnarray*}
||T_u|\varphi -P_t\varphi|^2||_\infty=||T_{u-t^2}T_{t^2}|\varphi
-P_t\varphi|^2||_\infty\leq ||T_{t^2}|\varphi
-P_t\varphi|^2||_\infty
\end{eqnarray*}
For $u\leq t^2$, denote $n$ the biggest integer smaller than $\frac
t{\sqrt {u}}$. We have
\begin{eqnarray*}
&&||T_u|\varphi -P_t\varphi|^2||_\infty^{\frac12}\\
&=&||T_{u}|\varphi -P_{\sqrt u}\varphi|^2||_\infty^{\frac12}+||T_{u}|P_{\sqrt u}\varphi -P_{2\sqrt u}\varphi|^2||_\infty^{\frac12}\\
&&+\cdot\cdot\cdot||T_{u}|P_{n\sqrt u}\varphi -P_{(n-1)\sqrt u}\varphi|^2||_\infty^{\frac12}+||T_{u}|P_t\varphi -P_{n\sqrt u}\varphi|^2||_\infty^{\frac12}\\
&\leq&||T_{u}|\varphi -P_{\sqrt u}\varphi|^2||_\infty^{\frac12}+||P_{\sqrt u}T_{u}|\varphi -P_{\sqrt u}\varphi|^2||_\infty^{\frac12}\\
&&+\cdot\cdot\cdot||P_{{n-1}\sqrt u}T_{u}|\varphi -P_{\sqrt u}\varphi|^2||_\infty^{\frac12}+||P_{n\sqrt u}T_{u-(t-n\sqrt u)^2}T_{(t-n\sqrt u)^2}|P_{t-n\sqrt u}\varphi -\varphi|^2||_\infty^{\frac12}\\
&\leq&2\frac t{\sqrt {u}}||\sup_tT_{t^2}|\varphi
-P_{t}\varphi|^2||_\infty^{\frac12}.
\end{eqnarray*}
Therefore,
\begin{eqnarray*}
 &&||P_{t}|\varphi -P_t\varphi
|^2||_\infty^{\frac12}\\
&\leq&(\frac 2{\sqrt{\pi }}\int_0^{t^2}
te^{-\frac{t^2}{4u}}u^{-\frac
32}\frac {t^2}udu)^\frac12||\sup_tT_{t^2}|\varphi -P_t\varphi|^2||_\infty^\frac12\\
&&+(\frac 1{2\sqrt{\pi }}\int_{t^2}^\infty
te^{-\frac{t^2}{4u}}u^{-\frac
32}du)^\frac12||\sup_tT_{t^2}|\varphi -P_t\varphi|^2||_\infty ^{\frac12}\\
&\leq&c||\sup_tT_{t^2}|\varphi -P_t\varphi|^2||_\infty
^{\frac12}.\qed
\end{eqnarray*}

\begin{proposition}\label{lemf}
Assume the underlying semigroup $(T_y)_y$ is quasi-monotone. Then,
for $\varphi\in L^2({\cal M})$, we have
\begin{eqnarray}
 ||\varphi
||_{BMO_c(P)}\approx\sup_{t,f}|\tau [\varphi^* (f-P_tf)]|,
\label{bmof}
\end{eqnarray}
where the supremum is taken for all $t>0$ and $f=bT_{t^2}^{\frac
12}(a)$ with $a,b\geq 0,\tau a\leq 1, \tau b^2\leq 1$.
%For
%$\varphi\in BMO_c(P)$ with a bounded sequence $\varphi_{\alpha}\in
%BMO_c(P)\cap L^2({\cal M})$ P-converges to $\varphi$,
%$\lim_{\alpha}\tau [\varphi_{\alpha}^*(f-P_tf)]$  exists and we
%understand $\tau [\varphi^* (f-P_tf)]$ as $\lim_{\alpha}\tau
%[\varphi^*(f-P_tf)]$.
\end{proposition}

{\bf Proof.}  Fix $t,\varphi\in L^2({\cal M}) ,$
\begin{eqnarray*}
||T_{t^2}|\varphi -P_t\varphi |^2||_\infty  &=&\sup_{a\geq 0,\tau
a\leq
1}\tau (aT_{t^2}|\varphi -P_t\varphi |^2) \\
&=&\sup_{a\geq 0,\tau a\leq 1}\tau (T_{t^2}(a)|\varphi -P_t\varphi |^2) \\
&=&\sup_{a\geq 0,\tau a\leq 1}\tau |(\varphi -P_t\varphi
)(T_{t^2}(a))^{\frac 12}|^2 \\
&=&\sup_{a\geq 0,\tau a\leq 1}\sup_{b\geq 0,\tau b^2\leq 1}(\tau
[b(T_{t^2}(a))^{\frac 12}(\varphi^* -P_t\varphi^* )])^2\\
&=&\sup_{a\geq 0,\tau a\leq 1}\sup_{b\geq 0,\tau b^2\leq 1}(\tau [
(b(T_{t^2}(a))^{\frac 12}-P_t[b(T_{t^2}(a))^{\frac
12}])\varphi^*])^2.
\end{eqnarray*}
Let
\begin{eqnarray}
f=b(T_{t^2}(a))^{\frac 12}.  \label{fg}
\end{eqnarray}
Then we get
\[
||\varphi ||_{BMO_c(P)}\approx \sup_t||T_{t^2}|\varphi -P_t\varphi
|^2||_\infty ^{\frac 12}=\sup_t\sup_{a\geq 0,\tau a\leq
1}\sup_{b\geq 0,\tau b^2\leq 1}\tau [\varphi^* (f-P_tf)].\qed
\]
%From the argument above, we easily see $\lim_{\alpha}\tau
%[\varphi_{\alpha}^*(f-P_tf)]$ exists if $\varphi_{\alpha}$
%P-converges and
%$$||w^*\lim_{\alpha} T_{t^2}|\varphi_{\alpha} -P_t\varphi_{\alpha}
%|^2||_\infty\approx\sup_{t,f}\lim_{\alpha}\tau [\varphi_{\alpha}^*
%(f-P_tf)].$$ The weak-$*$ limit of $T_{t^2}|\varphi_{\alpha}
%-P_t\varphi_{\alpha} |^2$ exists because of $(\ref{Py>Ty2}).$\qed

We will show $f-P_yf$ is in $H_c^1(P)$ with norm smaller than $c$.

\begin{proposition}
Given $t>0$, let
\begin{eqnarray*}
P_s^a &=&\int_0^{t^2}se^{-\frac{s^2}{4u}}u^{-\frac 32}T_udu; \\
P_s^b &=&\int_{t^2}^\infty se^{-\frac{s^2}{4u}}u^{-\frac 32}T_udu.
\end{eqnarray*}
Then, for any $0<s<\infty$, we have
\begin{eqnarray} P_s^b\leq c\frac stP_t^b  \label{psb}
\end{eqnarray}
and
\begin{eqnarray}
T_{t^2}P_s^a\leq 2^{\alpha} T_{t^2}, \label{psb1}
\end{eqnarray}
for ($T_y)_y$ quasi-decreasing with index ${\alpha}$;
\begin{eqnarray}
T_{t^2}P_s^a\leq 2^{\alpha} T_{2t^2}, \label{psb2}
\end{eqnarray}
for ($T_y)_y$ quasi-increasing with index ${\alpha}$.
\end{proposition}

{\bf Proof.} (\ref{psb}) is easy to verify by the facts that
$e^{-\frac{s^2}{4u}}$ decreases with respect to $s$ and
$e^{-\frac{s^2}{4u}}\approx e^{-\frac{t^2}{4u}}$ for any
$u>t^2,s<t$. Note the quasi decreasing (increasing) property implies
$T_{t^2+u}\leq 2^{\alpha} T_{t^2}$ ($T_{t^2+u}\leq 2^{\alpha}
T_{2t^2}$) for all $u<t^2$ respectively.  (\ref{psb1}) and
(\ref{psb2}) follow by the inequality
$\int_0^{t^2}se^{-\frac{s^2}{4u}}u^{-\frac 32}du\leq 1$. \qed

\begin{proposition}\label{dTyTy}
For ($T_y)_y=e^{yL}$ quasi-decreasing, we have
\begin{eqnarray}
-c_{\alpha} \frac{T_{\frac{2y}3}}y\leq \frac{\partial T_y}{\partial
y}\leq {\alpha} \frac{T_y}y. \label{LTyd} \end{eqnarray} For
($T_y)_y$ quasi-increasing, we have
\begin{eqnarray}
-{\alpha} \frac{T_y}y\leq \frac{\partial T_y}{\partial y}\leq
c_{\alpha} \frac{T_{2y}}y. \label{LTyi}
\end{eqnarray}
\end{proposition}

{\bf Proof. }Assume $\frac{T_y}{y^{\alpha} }$ decreasing, taking
derivative with respect to $y,$ we get
\[
\frac{\partial T_y}{\partial y}- {\alpha} \frac{T_y}y\leq0.
\]
which is the second inequality of (\ref{LTyd}). By using it, we get
\[
(-\frac{\partial T_y}{\partial y}+3{\alpha}
\frac{T_y}y)=(-\frac{\partial T_{\frac y3}}{\partial \frac y3}
+{\alpha} \frac{T_{\frac y3}}{\frac y3})T_{\frac{2y}3}\leq
2^{\alpha} (-\frac{\partial T_{\frac y3}}{\partial \frac
y3}+3{\alpha} \frac{T_{\frac y3}}y)T_s,
\]
for $\frac y3\leq s\leq \frac{2y}3.$ Taking integral for $s$ from
$\frac y3$ to $\frac{2y}3,$ we get
\begin{eqnarray*}
\frac y3(-\frac{\partial T_{y}}{\partial y}+3{\alpha} \frac{T_y}y) &\leq &\int_{\frac y3}^{\frac{2y}%
3}2^{\alpha} (-\frac{\partial T_{\frac y3}}{\partial \frac y3}+3{\alpha} \frac{T_{\frac y3}}y)T_sds \\
%&=&\int_{\frac y3}^{\frac{2y}3}2^{\alpha} (-LT_{\frac y3+s}+3{\alpha} \frac{%
%T_{\frac y3+s}}y)ds \\
&=&2^{\alpha}(\int_{\frac y3}^{\frac{2y}3}-\frac{\partial T_{\frac
y3+s}}{\partial s}
ds+\int_{\frac y3}^{\frac{2y}3}3{\alpha} \frac{T_{\frac y3+s}}yds) \\
&\leq &2^{\alpha}(-T_y+T_{\frac{2y}3}+\int_{\frac y3}^{\frac{2y}3}
3{\alpha} (\frac 32)^{\alpha}\frac{
T_{\frac{2y}3}}yds) \\
&=&-2^{\alpha} T_y+(3^{{\alpha}} {\alpha}
+2^{\alpha})T_{\frac{2y}3}.
\end{eqnarray*}
Therefore
\[
\frac{\partial T_{y}}{\partial y}\geq \frac{3(2^{\alpha}
+{\alpha})T_y-3(3^{{\alpha}} {\alpha}
+2^{\alpha})T_{\frac{2y}3}}y\geq -\frac{3(3^{\alpha}{\alpha}
+2^{\alpha}) T_{\frac{2y}3}}y.
\]
That is the first inequality of (\ref{LTyd}). The proof for
quasi-increasing semigroup is similar.\qed

\begin{lemma}\label{0tg}
Assume $(T_t)_t$ is quasi-monotone and satisfies the $L^{\frac12}$
condition (\ref{L1/2}). Then, for any $t>0$ and $f$ given as in
(\ref{fg}),
\begin{eqnarray}
\tau (\int_0^tT_{t^2}|\frac{\partial P_sf}{\partial s}|^2sds)^{\frac
12}\leq c_{\alpha} .  \label{0t}
\end{eqnarray}
\end{lemma}

{\bf Proof.} We only prove (\ref{0t}) for quasi-decreasing
semigroups. The proof for quasi-increasing ones is similar and
easier. For any positive element $x$ in $L^\infty({\cal M}),$ by
(\ref{Gamma>partial}) and (\ref{carrL}),
\begin{eqnarray*}
&&\tau (x\int_0^tT_{t^2}|\frac{\partial P_sf}{\partial s}|^2sds) \\
&=&\tau \int_0^tT_{t^2}(x)|\frac{\partial P_sf}{\partial s}|^2sds \\
&\leq &\tau \int_0^tT_{t^2}(x)\widetilde{\Gamma }(P_sf,P_sf)sds \\
&=&\tau \int_0^tT_{t^2}(x)(L+\frac{\partial ^2}{\partial
s^2})|P_sf|^2sds
\\
&=&\tau \int_0^tLT_{t^2}(x)|P_sf|^2sds+\tau[
T_{t^2}(x)\int_0^t\frac{\partial
^2}{\partial s^2}|P_sf|^2sds] \\
&=&I+II.
\end{eqnarray*}
For II, using of ``integration by parts",
\begin{eqnarray*}
II &=&\left. \tau T_{t^2}(x)s\frac \partial {\partial
s}|P_sf|^2\right| _{s=t}-\left. \tau T_{t^2}(x)s\frac \partial
{\partial s}|P_sf|^2\right|
_{s=0}-\tau \int_0^tT_{t^2}(x)\frac \partial {\partial s}|P_sf|^2ds \\
&=&\left. \tau T_{t^2}(x)s(P_sf\frac \partial {\partial
s}P_sf+(\frac
\partial {\partial s}P_sf)P_sf)\right|_{s=t}-0-\tau
T_{t^2}(x)\int_0^t\frac \partial {\partial s}|P_sf|^2ds \\
&\leq &\tau \left. T_{t^2}(x)t(\frac 1t|P_tf|^2+t|\frac \partial
{\partial s}P_sf|^2\right| _{s=t})-\tau
T_{t^2}(x)(|P_tf|^2-|f|^2)\\
&= &\tau \left. T_{t^2}(x)t^2|\frac \partial {\partial
s}P_sf|^2\right| _{s=t}+\tau T_{t^2}(x)|f|^2.
\end{eqnarray*} By the identity (\ref{idpy}), we get
\[
\left.\frac \partial {\partial s}P_sf\right| _{s=t}=\int_0^\infty (1-\frac{t^2}{2u})e^{-\frac{%
t^2}{4u}}u^{-\frac 32}T_ufdu.
\]
Then \begin{eqnarray*} \left.t^2|\frac \partial {\partial
s}P_sf|^2\right| _{s=t}&\leq& 2|\int_0^{2t^2}t(1-\frac{t^2}{2u})e^{-\frac{t^2}{4u}%
}u^{-\frac 32}T_ufdu|^2+2|\int_{2t^2}^\infty t(1-\frac{t^2}{2u})e^{\frac{t^2%
}{4u}}u^{-\frac 32}T_ufdu|^2\\
%&\leq& 4|\int_0^{2t^2}t(1+\frac{t^2}{2u})e^{-\frac{t^2}{4u}
%}u^{-\frac 32}T_ufdu|^2+2|\int_{2t^2}^\infty t(1-\frac{t^2}{2u})e^{\frac{t^2
%}{4u}}u^{-\frac 32}T_ufdu|^2\\
&\leq& c\int_0^{2t^2}t(1+\frac{t^2}{2u})e^{-\frac{t^2}{4u}
}u^{-\frac 32}|T_uf|^2du+2|\int_{2t^2}^\infty t(1-\frac{t^2}{2u})e^{\frac{t^2%
}{4u}}u^{-\frac 32}T_ufdu|^2\\
&\leq& c\int_0^{2t^2}t(1+\frac{t^2}{2u})e^{-\frac{t^2}{4u}
}u^{-\frac 32}T_u|f|^2du+2|\int_{2t^2}^\infty
t(1-\frac{t^2}{2u})e^{\frac{t^2 }{4u}}u^{-\frac 32}T_ufdu|^2.
\end{eqnarray*}
Therefore,
\begin{eqnarray*}
&&II\\
&\leq &\tau [xT_{t^2}(c\int_0^{2t^2}t(1+\frac{t^2}{2u})e^{-\frac{t^2}{4u}%
}u^{-\frac 32}T_u|f|^2du+2|\int_{2t^2}^\infty t(1-\frac{t^2}{2u})e^{-\frac{t^2%
}{4u}}u^{-\frac 32}T_ufdu|^2+|f|^2)] \\
&=&\tau [x(c\int_0^{2t^2}t(1+\frac{t^2}{2u})e^{-\frac{t^2}{4u}%
}u^{-\frac 32}T_{t^2+u}|f|^2du+T_{t^2}|f|^2+2T_{t^2}|\int_{2t^2}^\infty t(1-\frac{t^2}{2u})e^{-\frac{t^2%
}{4u}}u^{-\frac 32}T_ufdu|^2)] \\
&\leq&\tau [x(c\int_0^{2t^2}t(1+\frac{t^2}{2u})e^{-\frac{t^2}{4u}%
}u^{-\frac 32}3^{\alpha} T_{t^2}|f|^2du+T_{t^2}|f|^2+2T_{t^2}|\int_{2t^2}^\infty t(1-\frac{t^2}{2u})e^{-\frac{t^2%
}{4u}}u^{-\frac 32}T_ufdu|^2)] \\
&\leq &\tau [x(c_{\alpha} T_{t^2}|f|^2+2
T_{t^2}|T_{t^2}\int_{2t^2}^\infty t(1-\frac{t^2}{2u})e^{-\frac{t^2}{4u}%
}u^{-\frac 32}T_{u-t^2}fdu|^2)].
\end{eqnarray*}
Set
\[
h=\int_{2t^2}^\infty t(1-\frac{t^2}{2u})e^{-\frac{t^2}{4u}}u^{-\frac
32}T_{u-t^2}fdu.
\]
We get
\begin{eqnarray}
||h||_{L^1}&\leq& c||f||_{L^1}\leq c  \label{h} \\
II&\leq& \tau [x(c_{{{\alpha}}} T_{t^2}|f|^2+2
T_{t^2}|T_{t^2}h|^2)]. \nonumber
\end{eqnarray}
For I, by (\ref{LTyd}) and (\ref{psb1}), we have
\begin{eqnarray*}
I&=&\tau \int_0^t\left.\frac{\partial T_{y}}{\partial y}\right|_{y=t^2}(x)|P_sf|^2sds \\
&\leq &\tau \int_0^tc_{\alpha} \frac{T_{t^2}}{t^2}(x)|P_sf|^2sds \\
&\leq &2c_{\alpha} \tau
\int_0^t\frac{T_{t^2}}{t^2}(x)|P_s^af|^2sds+2c_{\alpha}
\tau \int_0^t\frac{T_{t^2}}{t^2}(x)|P_s^bf|^2sds \\
&\leq &2c_{\alpha} \tau
\int_0^tP_s^a\frac{T_{t^2}}{t^2}(x)|f|^2sds+2c_{\alpha}
\tau [\frac{T_{t^2}}{t^2}(x)\int_0^t|P_s^bf|^2sds] \\
&\leq &c_{\alpha} 2^{\alpha} \tau[T_{t^2}(x)|f|^2]+2c_{\alpha} \tau [\frac{T_{t^2}}{t^2}%
(x)\int_0^t|P_s^bf|^2sds].
\end{eqnarray*}
By (\ref{psb}), we have ($P_s^bf)^{\frac 12}=c(P_t^bf)^{\frac
12}u_s$ for $s<t$ with some partial contraction $u_s.$ Then
\begin{eqnarray*}
I &\leq &c_{\alpha} 2^{\alpha} \tau [T_{t^2}(x)|f|^2]+2c_{\alpha} \tau [\frac{T_{t^2}}{t^2%
}(x)\int_0^t(P_s^bf)^{\frac 12}(P_s^bf)(P_s^bf)^{\frac 12}sds] \\
&=&c_{\alpha} 2^{\alpha} \tau [x(T_{t^2}|f|^2)]+2c_{\alpha} \tau [\frac{T_{t^2}}{t^2}%
(x)\int_0^t(P_t^bf)^{\frac 12}u_s(P_s^bf)u_s^*(P_t^bf)^{\frac 12}sds] \\
&=&c_{\alpha} 2^{\alpha} \tau[x(T_{t^2}|f|^2)]+2c_{\alpha} \tau [x(\frac{T_{t^2}}{t^2}%
(P_t^bf)^{\frac 12}\int_0^tu_s(P_s^bf)u_s^*sds(P_t^bf)^{\frac 12})].
\end{eqnarray*}
Let
\[
g=\int_0^tu_s(P_s^bf)u_s^*sds.
\]
We see
\[
||g||_{L^1}\leq \frac{t^2}2||f||_{L^1}.
\]
Combining the estimations for I and II, we get
\begin{eqnarray*}
&&\tau x\int_0^tT_{t^2}|\frac{\partial P_sf}{\partial s}|^2sds \\
&\leq &c_{\alpha} \tau [x(T_{t^2}|f|^2+ T_{t^2}|T_{t^2}h|^2+
T_{t^2}[(P_t^bf)^{\frac 12}\frac{g}{t^2}(P_t^bf)^{\frac 12}])].
\end{eqnarray*}
By the arbitrariness of $x,\,$we get
\begin{eqnarray*}
\int_0^tT_{t^2}|\frac{\partial P_sf}{\partial s}|^2sds &\leq
&c_{\alpha}
(T_{t^2}|f|^2+T_{t^2}|T_{t^2}h|^2+T_{t^2}[(P_t^bf)^{\frac 12}\frac
g{t^2}(P_t^bf)^{\frac12}].
\end{eqnarray*}
Note $f=bT_{t^2}^{\frac 12}(a)$ and $P_t^bf$ is in form of
$T_{t^2}z$ with $||z||_{L^1}\leq ||f||_{L^1}\leq1$. Using the $L^{\frac 12}$ assumption for $%
T_y,$ we get by (\ref{h}),
\[
\tau (\int_0^tT_{t^2}|\frac{\partial P_sf}{\partial s}|^2sds)^{\frac
12}\leq
c_{{{\alpha}}} (||a||_{L^1}||b^2||_{L^1}+||h||_{L^1}^2+||f||_{L^1}^2)\leq c_{{%
{\alpha}}} .\qed
\]

\begin{lemma}\label{ks2}
Assume that $(T_t)_t$ is a positive semigroup as in section 1.2
(1.4). Then, for $f\geq0, ||f||_{L^1}\leq1$, we have
\[
\tau (\int_t^\infty |T_{ks^2}\frac{\partial P_s}{\partial s}
(f-P_tf)|^2sds)^{\frac 12}\leq c_k.
\]
for any positive scalar $k,t$.
\end{lemma}

{\bf Proof.} Let
\[
Q_s= T_{ks^2}\frac{\partial P_s}{\partial s}(f-P_tf)
\]
The identity (\ref{idpy}) yields
\begin{eqnarray*}
Q_s(f)&=&T_{ks^2}\frac{\partial P_s}{\partial
s}f-T_{ks^2}\frac{\partial P_{s+t}}{\partial s}f
\\
&=&\int_0^\infty [(1-\frac{s^2}{2u})e^{-\frac{s^2}{4u}}-(1-\frac{(s+t)^2}{2u%
})e^{-\frac{(s+t)^2}{4u}}]u^{-\frac 32}T_{u+ks^2}fdu\\
&=&\int_0^\infty \psi_s(u)T_{u+ks^2}fdu
\end{eqnarray*}
with
\[
\psi_s(u)=[(1-\frac{s^2}{2u})e^{-\frac{s^2}{4u}}-(1-\frac{(s+t)^2}{2u})e^{-\frac{(s+t)^2}{4u}}]u^{-\frac
32}
\]
Since $s\geq t$, we have $\psi_s(u)\thickapprox t\frac
{\partial}{\partial s}
[(1-\frac{s^2}{2u})e^{-\frac{s^2}{4u}}]u^{-\frac 32}$ and
$$|\psi_s(u)|\leq c\frac ts u^{-\frac 32}e^{-\frac{s^2}{4u}}\leq c_k\frac ts (u+ks^2)^{-\frac 32}.$$
Let
\[
 R_t(f)=\int_{kt^2}^\infty tu^{-\frac 32}T_u(f)du
\]
Noting that $kt^2\leq ks^2$, we have $$-c_k\frac {R_t(f)}s\leq
Q_s(f)\leq c_k\frac {R_t(f)}s$$ Then there exist partial
contractions $u_s$ such that
\[
|Q_s(f)|^{\frac 12}=(c_k\frac {R_t(f)}s)^{\frac 12}u_s
\]
Then
%\begin{eqnarray*}
%&&|Q_s(f)|^2 \\
%&=&\frac {c_k}s[(R_{t}f)^{\frac 12}u_s|Q_s(f)|u_s^{*}(R_{t}f)^{\frac
%12}].
%\end{eqnarray*}
%And
\begin{eqnarray*}
&&\tau (\int_t^\infty |Q_s(f)|^2sds)^{\frac 12}\\
&=&c_k^\frac12\tau (\int_t^\infty [(R_{t}f)^{\frac
12}u_s|Q_s(f)|u_s^{*}(R_{t}f)^{\frac 12}]ds)^{\frac 12} \\
&=&c_k^\frac12\tau[(R_{t}f)^{\frac 12}\int_t^\infty
u_s|Q_s(f)|u_s^{*}ds(R_{t}f)^{\frac 12}]^{\frac 12}.
\end{eqnarray*}
Note that
\begin{eqnarray*}
 ||R_t(f)||_{L^1}&\leq& 2k^{-\frac12};\\
||\int_t^\infty u_s|Q_s(f)|u_s^{*}ds||_{L^1}&\leq& \int_t^\infty
||Q_s(f)||_{L^1}ds\\
&\leq&\int_t^\infty\int_0^\infty|\psi_s(u)|duds\\
&\leq&c_k\int_t^\infty\frac t{s^2}ds= c_k.\\
\end{eqnarray*}
By H\"{o}lder's inequality, we get
\begin{eqnarray*}\tau (\int_t^\infty |T_{ks^2}\frac{\partial
P_s}{\partial s} (f-P_tf)|^2sds)^{\frac 12}&\leq& c_k^{\frac
12}||R_t(f)||_{L^1}^{\frac 12}||\int_t^\infty
u_sQ_s(f)u_s^{*}ds||_{L^1}^{\frac 12}\\&\leq&
c_k.\qed\end{eqnarray*}

%Therefore,
%\begin{eqnarray*}
%||f||_{H_G^1} &=&\tau (\int_0^\infty |\frac{\partial P_s}{\partial s}%
%(f-P_tf)|^2sds)^{\frac 12} \\
%&=&\tau T_{t^2}(\int_0^\infty |\frac{\partial P_s}{\partial s}%
%(f-P_tf)|^2sds)^{\frac 12} \\
%&\leq &\tau (\int_0^\infty T_{t^2}|\frac{\partial P_s}{\partial s}%
%(f-P_tf)|^2sds)^{\frac 12} \\
%&\leq &2\tau (\int_0^tT_{t^2}|\frac{\partial P_s}{\partial s}%
%(f-P_tf)|^2sds)^{\frac 12}+2\tau (\int_t^\infty T_{t^2}|\frac{\partial P_s}{%
%\partial s}(f-P_tf)|^2sds)^{\frac 12} \\
%&\leq &c_{\alpha} .
%\end{eqnarray*}
%We finished the proof of Lemma \ref{ks2}.

\begin{lemma}\label{ks8}
Assume that $(T_y)_y$ is quasi monotone with index ${\alpha}$ and satisfy the $L^{\frac 12}$ condition (\ref{fg}%
). There exists a constant $k\leq 4$ depending only on ${\alpha}$
such that
\begin{eqnarray*}
\tau (\int_0^\infty T_{s^2}|\frac{\partial P_s}{\partial
s}g|^2sds)^{\frac 12}\leq c_{\alpha}\tau (\int_0^\infty |T_{\frac
{ks^2}8}\frac{\partial P_{s}}{
\partial s}g|^2sds)^{\frac 12},
\end{eqnarray*}
for any $g$.
\end{lemma}

{\bf Proof.} We will prove only for quasi-increasing $(T_y)_y$ since
the proof for quasi-decreasing ones is easier (and similar) for this
Lemma. By Proposition \ref{alpha}, we can find a constant
$c_{\alpha}\geq1$ such that
\begin{eqnarray}
\tau (\int_0^\infty T_{2s^2}|A_s|^2\frac{ds}s)^{\frac 12}\leq
c_{{\alpha} }\tau (\int_0^\infty T_{(\frac s2)^2}|A_s|^2\frac{ds}
s)^{\frac 12}.\label{calpha}
\end{eqnarray}
Choose scalar $k\leq4$ such that
$$
c^2_{{\alpha}}2^{\alpha}\int_0^{ks^2}se^{-\frac{s^2}{4u}}u^{-\frac
32} ds\leq \frac1{16}
$$
Set
\begin{eqnarray*}
P_{s}^c &=&\int_0^{ks^2}se^{-\frac{s^2}{4u}}u^{-\frac 32}T_udu; \\
P_{s}^d &=&\int_{ks^2}^\infty se^{-\frac{s^2}{4u}}u^{-\frac
32}T_udu.
\end{eqnarray*}
Then, for $(T_s)_s$ quasi-increasing,
\begin{eqnarray}
T_{4s^2}P_{s}^c &=&\int_0^{ks^2}se^{-\frac{s^2}{4u}}u^{-\frac 32}T_{u+4s^2}du \nonumber\\
&\leq&\int_0^{ks^2}se^{-\frac{s^2}{4u}}u^{-\frac 32}du2^{\alpha} T_{8s^2} \nonumber\\
&\leq&\frac1{16c^2_{{\alpha}}}T_{8s^2}.\label{tsps}
\end{eqnarray}
By (\ref{tsps}), for $t$ fixed, we get
\begin{eqnarray*}
&&\tau (\int_0^\infty T_{s^2}|\frac{\partial P_s}{\partial
s}g|^2sds)^{\frac
12} \\
&=&\tau (\int_0^\infty T_{s^2}|P_{\frac s2}\frac{\partial P_{\frac
s2}}{
\partial s}g|^2sds)^{\frac 12} \\
&\leq &\tau (\int_0^\infty T_{s^2}|P_{\frac s2}^c\frac{\partial P_{\frac s2}%
}{\partial \frac s2}g|^2sds)^{\frac 12}+\tau (\int_0^\infty T_{s^2}|P_{\frac s2}^d%
\frac{\partial P_{\frac s2}}{\partial \frac s2}g|^2sds)^{\frac 12} \\
&\leq &\tau (\int_0^\infty T_{s^2}P_{\frac s2}^c|\frac{\partial P_{\frac s2}%
}{\partial \frac s2}g|^2sds)^{\frac 12}+\tau (\int_0^\infty T_{s^2}|P_{\frac s2}^d%
\frac{\partial P_{\frac s2}}{\partial \frac s2}g|^2sds)^{\frac 12} \\
&\leq&\tau (\int_0^\infty
\frac1{16c_{\alpha}^2}T_{2s^2}|\frac{\partial P_{\frac s2}}{\partial
\frac s2}g|^2sds)^{\frac 12}+\tau (\int_0^\infty T_{s^2}|P_{\frac
s2}^d\frac{\partial P_{\frac s2}}{\partial \frac s2}g|^2sds)^{\frac
12}
\end{eqnarray*}
Applying  (\ref{calpha}), we get
\begin{eqnarray*}
\tau (\int_0^\infty T_{s^2}|\frac{\partial P_s}{\partial
s}g|^2sds)^{\frac 12}&\leq &\frac14\tau (\int_0^\infty T_{(\frac
s2)^2}|\frac{\partial P_{\frac s2}}{\partial \frac
s2}g|^2sds)^{\frac 12}+\tau (\int_0^\infty T_{s^2}|P_{\frac
s2}^d\frac{\partial P_{\frac s2}}{\partial \frac s2}g|^2sds)^{\frac 12} \\
&= &\frac 12\tau (\int_0^\infty T_{s^2}|\frac{\partial P_s}{\partial s}%
f|^2sds)^{\frac 12}+\tau (\int_0^\infty T_{s^2}|P_{\frac s2}^d\frac{%
\partial P_{\frac s2}}{\partial \frac s2}g|^2sds)^{\frac 12}.
\end{eqnarray*}
Thus
\begin{eqnarray}
\tau (\int_0^\infty T_{s^2}|\frac{\partial P_s}{\partial
s}g|^2sds)^{\frac 12}\leq 2\tau (\int_0^\infty T_{s^2}|P_{\frac
s2}^d\frac{\partial P_{\frac s2}}{\partial \frac s2}g|^2sds)^{\frac
12}.\label{psd}
\end{eqnarray}
Let
\[
P_{s}^e=\int_{ks^2}^\infty se^{-\frac{s^2}{4u}}u^{-\frac 32}T_{u-
\frac{k^2s^2}2}du.
\]
Then
\begin{eqnarray*}
P_{s}^d&=&P_s^eT_{\frac {ks^2}2}.
\end{eqnarray*}
For $(T_s)_s$ quasi-increasing, we have
\begin{eqnarray*}
T_{t^2}P_{\frac s2}^e&=&\int_{\frac {ks^2}4}^{t^2} \frac
s2e^{-\frac{s^2}{16u}}u^{-\frac 32}T_{u-
\frac{k^2s^2}8+t^2}du+T_{t^2}\int_{t^2}^\infty \frac
s2e^{-\frac{s^2}{16u}}u^{-\frac 32}T_{u-
\frac{k^2s^2}8}du\\
&\leq&\int_{\frac {ks^2}4}^{t^2} \frac
s2e^{-\frac{s^2}{16u}}u^{-\frac 32}2^{\alpha}
T_{2t^2}du+T_{t^2}\int_{t^2}^\infty \frac
s2e^{-\frac{s^2}{16u}}u^{-\frac 32}T_{u-
\frac{k^2s^2}8}du\\
&\leq& 2^{\alpha}(T_{2t^2}+T_{t^2}\int_{t^2}^\infty \frac t2u^{-\frac 32}T_{u}du)\\
&\leq& 2^{\alpha}(T_{2t^2}+T_{t^2}P_{t}).
\end{eqnarray*}
for any $s\leq t$. Applying this inequality,  we have, for any
$(B_{s})_s$,
\begin{eqnarray*}
||(P_{\frac s2}^eB_{s})_s||_{{\cal T}_\infty^{(T_{s^2})}}&=&||\sup_tT_{t^2}\int_0^{t}|P_{\frac s2}^eB_{s}|^2\frac {ds}s||_\infty^{\frac12}\\
&\leq&||\sup_t\int_0^{t}T_{t^2}P_{\frac s2}^e|B_{s}|^2\frac {ds}s||_\infty^{\frac12}\\
&\leq&2^{\frac\alpha2}||\sup_tT_{2 t^2}\int_0^{t}|B_{s}|^2\frac
{ds}s||_\infty^\frac12
+2^{\frac\alpha2}||\sup_tP_{t}T_{t^2}\int_0^{t}|B_{s}|^2\frac {ds}s||_\infty^{\frac12}\\
&\leq&2^{\frac\alpha2}\sup_t||T_{t^2}\int_0^{t}|B_{s}|^2\frac {ds}s||_\infty^{\frac12}+2^{\frac\alpha2}\sup_t||T_{t^2}\int_0^{t}|B_{s}|^2\frac {ds}s||_\infty^{\frac12}\\
&\leq&c_{\alpha}||(B_{s})_s||_{{\cal T}_\infty^{(T_{s^2})}}.
\end{eqnarray*}
By the duality between tent spaces ${\cal T}_1^{(T_{s^2})}$ and
${\cal T}_\infty^{(T_{s^2})}$, which is implied by the assumption of
the lemma and Remark \ref{change}, we get
\begin{eqnarray}
||(P_{\frac s2}^eA_{s})_s||_{{\cal
T}_1^{(T_{s^2})}}&\leq&c_{\alpha}||(A_{s})_s||_{{\cal
T}_1^{(T_{s^2})}}. \label{pset1}
\end{eqnarray}
Applying (\ref{pset1}) to (\ref{psd}) and using Proposition
\ref{alpha}, we get
\begin{eqnarray}
\tau (\int_0^\infty T_{s^2}|\frac{\partial P_s}{\partial
s}g|^2sds)^{\frac 12}&\leq& 2\tau (\int_0^\infty T_{s^2}|P_{\frac
s2}^eT_{\frac {ks^2}8}\frac{\partial P_{\frac
s2}}{\partial \frac s2}g|^2sds)^{\frac 12}\nonumber\\
&\leq& 2c_{\alpha}\tau (\int_0^\infty T_{s^2}|T_{\frac
{ks^2}8}\frac{\partial P_{\frac s2}}{\partial \frac
s2}g|^2sds)^{\frac 12}\\
& \leq& c_{{\alpha}}\tau (\int_0^\infty T_{\frac
{ks^2}{32}}|T_{\frac {ks^2}8}\frac{\partial P_{\frac s2}}{\partial
\frac s2}g|^2sds)^{\frac 12}.\label{psd1}
\end{eqnarray}
By Lemma \ref{T2s}, there exists another constant $c_{\alpha}$ such
that
\begin{eqnarray}
&&\tau (\int_0^\infty T_{\frac {ks^2}{32}}|T_{\frac
{ks^2}8}\frac{\partial P_{\frac s2}}{
\partial \frac s2}g|^2sds)^{\frac 12}\nonumber\\
&\leq&c_{\alpha}[\tau (\int_0^\infty T_{\frac
{ks^2}{32}}|\frac{\partial P_{\frac s2}}{
\partial \frac s2}g|^2sds)^{\frac 12}]^{\frac 12}[\tau (\int_0^\infty |T_{\frac {ks^2}{32}}\frac{\partial P_{\frac s2}}{
\partial \frac s2}g|^2sds)^{\frac 12}]^{\frac 12}.\label{ct2s}
\end{eqnarray}
Combining (\ref{psd1}), (\ref{ct2s}) and applying Proposition \ref
{alpha} again, we get
\begin{eqnarray*}
\tau (\int_0^\infty T_{s^2}|\frac{\partial P_s}{\partial
s}g|^2sds)^{\frac 12}&\leq&c_{\alpha}[\tau (\int_0^\infty T_{\frac
{ks^2}{32}}|\frac{\partial P_{\frac s2}}{
\partial \frac s2}g|^2sds)^{\frac 12}]^{\frac 12}[\tau (\int_0^\infty |T_{\frac {ks^2}{32}}\frac{\partial P_{\frac s2}}{
\partial \frac s2}g|^2sds)^{\frac 12}]^{\frac 12}\\
&\leq&c_{\alpha}[\tau (\int_0^\infty T_{\frac
{s^2}{4}}|\frac{\partial P_{\frac s2}}{
\partial \frac s2}g|^2sds)^{\frac 12}]^{\frac 12}[\tau (\int_0^\infty |T_{\frac {ks^2}{32}}\frac{\partial P_{\frac s2}}{
\partial \frac s2}g|^2sds)^{\frac 12}]^{\frac 12}\\
&=&c_{\alpha}[\tau (\int_0^\infty T_{{s^2}}|\frac{\partial P_{s}}{
\partial s}g|^2sds)^{\frac 12}]^{\frac 12}[\tau (\int_0^\infty |T_{\frac {ks^2}8}\frac{\partial P_{ s}}{
\partial s}g|^2sds)^{\frac 12}]^{\frac 12}
\end{eqnarray*}
Therefore,
\begin{eqnarray*}
\tau (\int_0^\infty T_{s^2}|\frac{\partial P_s}{\partial
s}f|^2sds)^{\frac 12}\leq c_{\alpha}\tau (\int_0^\infty |T_{\frac
{ks^2}8}\frac{\partial P_{s}}{
\partial s}f|^2sds)^{\frac 12}.\qed
\end{eqnarray*}

\begin{theorem}
Assume that the underlying semigroup $(T_y)_y$ is quasi-monotone and
satisfies the $L^{\frac 12}$ condition (\ref{L1/2}). Then BMO$_c(P)
=(H_c^1(P))^{*}.$
\end{theorem}

{\bf Proof. } The relation BMO$_c(P) \subset (H_c^1(P))^{*}$ is
Theorem \ref{3.5}. We only need to show \begin{eqnarray} ||\varphi
||_{BMO_c(P)}\leq c||\varphi ||_{(H_c^1)^{*}},
\label{onlyneed}\end{eqnarray} for $\varphi\in L^2({\cal M})\cap
BMO_c(P)$. Once this is proved, by the proof of Theorem \ref{non1.2}
and the Hahn-Banach theory, any linear functional $\ell$ on
$H_c^1(P)$ is given by
\begin{eqnarray}
\ell(f)=\lim_{k}\tau \int_{0}^\infty s\frac{\partial P_s f}{\partial
s} (g_s^k)^* \frac{ds}s=\lim_{k,n}\tau f\int_{\frac 1n}^n
\frac{\partial P_s}{\partial s}(g_s^k)^* {ds}
\end{eqnarray}
for $f\in L^2({\cal M})\cap H_c^1(P)$ with $g^k\in {\cal
T}_\infty^{(T_{s^2})}\cap L^2({\cal M}, L_c^2)$ such that
$||(g^k_s)_s||_{{\cal T}_\infty}\leq c||\ell||.$

Let $$\varphi_{k,n}=\int_{\frac 1n}^n \frac{\partial P_s}{\partial
s}(g_s^{k})^* {ds}\in L^2({\cal M}).$$ Because of (\ref{onlyneed}),
we have
$$||\varphi_{k,n}||_{BMO_c(P)}\leq c||\int_{\frac 1n}^n \frac{\partial
P_s}{\partial s}(g_s^{k})^* {ds}||_{(H_c^1)^{*}}\leq
c||(g_s^{k})_s||_{{\cal T}_\infty^{(T_{s^2})}}\leq c||\ell||.$$
There exists a subsequence which P-converges to an element
$\varphi\in BMO_c(P)$ with
$$||\varphi||_{BMO_c(P)}\leq
\sup_{k,n}||\varphi_{k,n}||_{BMO_c(P)}\leq c||\ell||,$$ because the
unit ball of $L^\infty({\cal M})\otimes L^\infty({\Bbb R}_+)$ is
weak-$*$ compact.

We now prove (\ref{onlyneed}). Because of Proposition \ref{lemf}, we
only need to show
\[
g=f-P_tf\in H_c^1(P)
\]
for any $f$ given as in (\ref{fg}).

Let $k$ be the constant in Lemma \ref{ks8}, we have
\begin{eqnarray*}
||g||_{H_c^1(P)}&=&\tau(\int_0^\infty T_{s^2}|\frac {\partial P_sg}{\partial s}|^2sds)^{\frac12}\\
&\leq&c_{\alpha}\tau(\int_0^\infty|T_{\frac {ks^2}8}\frac{\partial P_sg}{\partial s}|^2sds)^{\frac12}\\
&\leq&c_{\alpha}\tau(\int_0^t|T_{\frac {ks^2}8}\frac {\partial P_sg}{\partial s}|^2sds)^{\frac12}+c_{\alpha}\tau(\int_t^\infty|T_{\frac {ks^2}8}\frac {\partial P_sg}{\partial s}|^2sds)^{\frac12}\\
&\leq&c_{\alpha}\tau(\int_0^tT_{\frac {ks^2}8}|\frac {\partial
P_sg}{\partial s}|^2sds)^{\frac12}
+c_{\alpha}\tau(\int_t^\infty|T_{\frac {ks^2}8}\frac {\partial
P_sg}{\partial s}|^2sds)^{\frac12}.
\end{eqnarray*}
From Lemma \ref{ks2}, we know the second term is smaller than $c_k$.

For the first term, if $(T_s)_s$ is quasi-increasing, since $k\leq
4$, we have
\begin{eqnarray*}
\tau(\int_0^tT_{\frac {ks^2}8}|\frac {\partial P_sg}{\partial
s}|^2sds)^{\frac12}
&=&\tau T_{\frac{t^2}2}(\int_0^tT_{\frac {ks^2}8}|\frac {\partial P_sg}{\partial s}|^2sds)^{\frac12}\\
&\leq&\tau(\int_0^tT_{\frac{t^2}2+\frac {ks^2}8}|\frac {\partial P_sg}{\partial s}|^2sds)^{\frac12}\\
&\leq&2^{\frac {\alpha}2}\tau(\int_0^tT_{t^2}|\frac {\partial
P_sg}{\partial s}|^2sds)^{\frac12}.
\end{eqnarray*}
 For quasi-decreasing $(T_s)_s$, we get similarly,
\begin{eqnarray*}
\tau(\int_0^tT_{\frac {ks^2}8}|\frac {\partial P_sg}{\partial
s}|^2sds)^{\frac12}
&=&\tau T_{{t^2}}(\int_0^tT_{\frac {ks^2}8}|\frac {\partial P_sg}{\partial s}|^2sds)^{\frac12}\\
&\leq&\tau(\int_0^tT_{{t^2}+\frac {ks^2}8}|\frac {\partial P_sg}{\partial s}|^2sds)^{\frac12}\\
&\leq&2^{\frac {\alpha}2}\tau(\int_0^tT_{t^2}|\frac {\partial
P_sg}{\partial s}|^2sds)^{\frac12}.
\end{eqnarray*}
Therefore,
\begin{eqnarray}
||g||_{H_c^1(P)}\leq c_{\alpha}\tau(\int_0^tT_{t^2}|\frac {\partial
P_sg}{\partial
s}|^2sds)^{\frac12}+c_{\alpha}\tau(\int_t^\infty|T_{\frac
{ks^2}8}\frac {\partial P_sg}{\partial s}|^2sds)^{\frac12}.
\label{gHs}
\end{eqnarray}
Applying Lemmas \ref{0tg} and \ref{ks2} to (\ref{gHs}), we get $
||g||_{H_c^1(P)}\leq c_{\alpha}.\qed $

Once again, if $(T_y)_y$ is  classical heat semigroup on
$\Bbb{R}^n$, $||f||_{H_c^1(P)}$ is equivalent to the classical Hardy space $H^1$ norm of $f$ and $%
||\varphi ||_{BMO_c(P)}$ is equivalent to the classical BMO norm of
$\varphi .$ We recover the duality between the classical $H^1$ and
BMO.

\section{$H^1, \mathrm{BMO}$ associated with general semigroups}

\setcounter{theorem}{0}\setcounter{equation}{0} In this section, we
discuss a pair of $H^1$,BMO-like spaces associated with general
semigroup $(T_s)_s$ satisfying the usual property (i)-(iv) listed in
Section 1.2. We do not assume that $(T_s)_s$ satisfies the
quasi-monotone conditions except in Theorem \ref{end}.

For $f\in L^2(\mathcal{M})$, let
\begin{eqnarray*}
S_T(f)&=&(\int_0^\infty T_s(|\frac {\partial T_s}{\partial
s}f|^2)sds)^{\frac12}, \\
G(f)&=&(\int_0^\infty |\frac {\partial T_s}{\partial s}f|^2sds)^{\frac12}, \\
C_t(f)&=&\int_0^tT_t|\frac {\partial T_s}{\partial s}f|^2sds.
\end{eqnarray*}

Set
\begin{eqnarray*}
||f||_{\mathcal{H}_{c,1}^{S}} &=&||S_T(f)||_{L^1},\\
||f||_{\mathrm{BMO}_{c}^{C} }&=&\sup_t||C_t(f)||_{L^\infty}^\frac12.
\end{eqnarray*}

Another $H^1$-norm associated with semigroups has been studied by
Stein ([St2]) in the commutative case and Junge, Le Merdy, Xu
([JLX])in the
noncommutative case. That is the norm defined for $f\in L^2(\mathcal{%
M})$ as
\[
||f||_{\mathcal{H}_{c,1}^{G}}=||G(f)||_{L_1},\ \ \ \forall 1\leq p<\infty .
\]
%We define the corresponding row spaces and column-row spaces in the usual
%way.
It is easy to see that
\begin{eqnarray}
||f||_{\mathcal{H}_{c,1}^{G}}\leq 2||f||_{\mathcal{H}_{c,1}^{S}}.
\label{HGS}
\end{eqnarray}
by (\ref{cp}).

\begin{theorem}
\label{end1.1} For any semigroup $(T_y)_{y\geq 0}$ satisfying
(i)-(iv) in Section 1.2, we have
\[
|{\tau }f\varphi ^{*}|\leq c||f||_{\mathcal{H}_{c,1}^S}||\varphi
|||_{BMO_c^C},
\]
for $f,\varphi\in L^2({\cal M})$.
\end{theorem}

We use the same idea of the proof of Theorem \ref{non1.1}. The
advantage of having specific elements allows us to make
modifications at some key points and remove the quasi-monotone
assumption for $(T_s)_s$. Set truncated square functions $S_s,G_s$
as follows:
\begin{eqnarray}
S_s &=&(\int_s^\infty T_{y-\frac s2}(|\frac {\partial T_{y+\frac
s2}}{\partial y}f|^2)ydy)^{\frac 12}  \label{Ss} \\
{G}_s &=&(\int_s^\infty |\frac {\partial T_{2y}}{2\partial
y}f|^2ydy)^{\frac 12}.  \label{Gs}
\end{eqnarray}
The square functions $S_s,G_s$ satisfy our key Lemma.

\begin{lemma}
\label{endlem1.1}
\begin{eqnarray}
{G}_s\leq S_s;  \label{endlemma1} \\
\frac{dT_s(S_s)}{ds}\geq 2T_{{\frac s2}}(\frac{dT_{\frac s2}(S_s)}{ds}),~%
\frac{dT_{\frac s2}(S_s)}{ds}\leq 0.  \label{endlemma}
\end{eqnarray}
\end{lemma}

\textbf{Proof}. (\ref{endlemma1}) is true because of the fact
\begin{eqnarray}
|\frac {\partial T_{2y}}{2\partial y}f|^2 =|T_{y-\frac s2}\frac
{\partial T_{y+\frac s2}}{\partial y}f|^2 \leq T_{y-\frac s2}(|\frac
{\partial T_{y+\frac s2}}{\partial y}f|^2)
\end{eqnarray}
for any $y\geq \frac s2$, which follows from (\ref{cp}).

By (\ref{cp}) again, we get $S_{s}\geq S_t$ for any $s\leq t$, then
\begin{eqnarray*}
T_{s+\Delta s}(S_{s+\Delta s})-T_{s}(S_s) &=&T_{{\frac s2}}[T_{\frac
{s+2\Delta s}2}(S_{s+\Delta s})-T_{\frac {s}2}(S_s)] \\
&\geq&T_{\frac s2}[T_{\frac {s+2\Delta s}2}(S_{s+2\Delta s})-T_{\frac
{s}2}(S_s)].
\end{eqnarray*}
Dividing by $\Delta s$ both the sides, we get the first inequality
of (\ref {endlemma}).

We go to prove the second inequality of (\ref{endlemma}). By (\ref{cp}) and (%
\ref{tys}), we get
\begin{eqnarray*}
&&T_{\frac{s+2\Delta s}2}S_{s+2\Delta s}-T_sS_s \\
&=&T_{\frac s2}T_{\Delta s}(\int_{s+2\Delta s}^\infty T_{y-\frac s2-\Delta
s}(|\frac {\partial T_{y+\frac s2+\Delta s}}{\partial y}f|^2)ydy)^{\frac
12}- T_{\frac s2}(\int_s^\infty T_{y-\frac s2}(|\frac {\partial T_{y+\frac
s2}}{\partial y}f|^2)ydy)^{\frac 12} \\
&\leq&T_{\frac s2}(\int_{s+2\Delta s}^\infty T_{y-\frac s2}(|\frac {\partial
T_{y+\frac s2+\Delta s}}{\partial y}f|^2)ydy)^{\frac 12}- T_{\frac
s2}(\int_s^\infty T_{y-\frac s2}(|\frac {\partial T_{y+\frac s2}}{\partial
y}f|^2)ydy)^{\frac 12} \\
&\leq&T_{\frac s2}(\int_{s+2\Delta s}^\infty T_{y+\frac{\Delta s}2-\frac
s2}(|\frac {\partial T_{y+\frac{\Delta s}2+\frac s2}}{\partial
y}f|^2)ydy)^{\frac 12}- T_{\frac s2}(\int_s^\infty T_{y-\frac s2}(|\frac
{\partial T_{y+\frac s2}}{\partial y}f|^2)ydy)^{\frac 12}.
\end{eqnarray*}
A change of variables implies that
\begin{eqnarray*}
&&T_{\frac s2}(\int_{s+2\Delta s}^\infty T_{y+\frac{\Delta s}2-\frac
s2}(|\frac {\partial T_{y+\frac{\Delta s}2+\frac s2}}{\partial
y}f|^2)ydy)^{\frac 12} \\
&=&T_{\frac s2}(\int_{s+2\Delta s+\frac{\Delta s}2}^\infty T_{u-\frac
s2}(|\frac {\partial T_{u+\frac s2}}{\partial y}f|^2)(u-\frac{\Delta s}%
2)du)^{\frac 12} \\
&\leq&T_{\frac s2}(\int_s^\infty T_{y-\frac s2}(|\frac {\partial T_{y+\frac
s2}}{\partial y}f|^2)ydy)^{\frac 12}.
\end{eqnarray*}
Then
\begin{eqnarray*}
T_{\frac{s+2\Delta s}2}S_{s+2\Delta s}-T_sS_s\leq 0.
\end{eqnarray*}
Taking $\Delta s\rightarrow 0$ we prove the second inequality of
(\ref {endlemma}).\qed

\begin{lemma}
\label{fphis} For any semigroup $(T_y)_{y\geq 0}$ satisfying
(i)-(iv) in Section 1.2, we have
\[
|{\tau }\int_0^\infty \frac{\partial T_{3s}f}{\partial s}\varphi
_s^{*}sds|\leq {3}\sup_y||T_{\frac y2}(\int_0^y|\varphi _s|^2)sds||_\infty
^{\frac 12}||G(f)||_1^{\frac 12}||S(f)||_1^{\frac 12}
\]
for and $f\in L^2({\cal M})$ and any family $(\varphi _s)_s\in {\cal
T}_\infty^{(T_{\frac y2})}$.
\end{lemma}

\textbf{Proof.} We can assume $G_s$ invertible by approximation. By (\ref{cp}%
), (\ref{endlem1.1}) and the Cauchy-Schwarz inequality, we get
\begin{eqnarray*}
|{\tau}\int_0^\infty \frac {\partial T_{3s}f}{\partial s}\varphi^*_ssds|
&=&3|{\tau}\int_0^\infty T_s\frac {\partial T_{2s}f}{2\partial s}\varphi^*_ssds| \\
&=&3|{\tau}\int_0^\infty \frac {\partial T_{2s}f}{2\partial s}T_{s}\varphi^*_ssds| \\
&\leq&3({\tau}\int_0^\infty |\frac {\partial T_{2s}f}{2\partial s}|^2{G}%
_s^{-1}sds )^{\frac 12}({\tau}\int_0^\infty |T_{s}\varphi_s|^2{G_s}%
sds)^{\frac 12} \\
&\leq&3({\tau}\int_0^\infty |\frac {\partial T_{2s}f}{2\partial s}|^2s{G}%
_s^{-1}ds )^{\frac 12}({\tau}\int_0^\infty |T_s\varphi_s|^2{S_s}sds)^{\frac
12} \\
&\stackrel{def}{=}&3I^{\frac 12}II^{\frac 12}.
\end{eqnarray*}

For $I,$ we have
\begin{eqnarray*}
I ={\tau}\int_0^\infty -\frac{\partial {G}_s^2}{\partial s}{G}_s^{-1}ds =2{%
\tau}\int_0^\infty -\frac{\partial {G}_s}{\partial s}ds =2||G_0||_1.
\end{eqnarray*}
For $II,$ by (\ref{cp}) and use the identity $T_{s}(S_s)=\int_s^\infty-\frac{%
\partial T_{y}(S_y)}{\partial y}dy$ we have
\begin{eqnarray}
II\leq{\tau}\int_0^\infty T_s|\varphi_s|^2S_ssds
&=&{\tau}\int_0^\infty |\varphi_s|^2T_{s}(S_s)sds  \nonumber \\
&=&{\tau}\int_0^\infty |\varphi_s|^2s\int_s^\infty-\frac{\partial T_{y}(S_y)%
}{\partial y}dyds  \nonumber \\
&=&-\tau\int_0^\infty \int_0^t|\varphi_s|^2sds\frac{\partial T_{y}(S_y)}{%
\partial y}dy.  \label{II}
\end{eqnarray}
Substituting (\ref{endlemma}) to (\ref{II}), we get
\begin{eqnarray*}
II&\leq&-2{\tau}\int_0^\infty \int_0^t|\varphi_s|^2sdsT_{\frac y2}(\frac{
\partial T_{\frac {y}2}(S_y)}{\partial y})dy \\
&=&-2{\tau}\int_0^\infty T_{\frac y2}\int_0^t|\varphi_s|^2sds\frac{ \partial
T_{\frac {y}2}(S_y)}{\partial y}dy \\
&\leq &2\sup_y||T_{\frac y2}(\int_0^y|\varphi_s|^2sds)||_\infty {\tau}
\int_0^\infty -\frac{\partial T_{\frac y2}(S_y)}{\partial y}dy \\
&=&2\sup_y||T_{\frac y2}(\int_0^y|\varphi_s|^2sds)||_\infty ||T_0(S_0)||_1 \\
&=&2\sup_y||T_{\frac y2}(\int_0^y|\varphi_s|^2sds)||_\infty ||S(f)||_1
\end{eqnarray*}
Combining the estimates of I and II, we get
\begin{eqnarray*}
|{\tau}\int_0^\infty \frac {\partial T_{3s}f}{\partial s}\varphi_ssds|\leq {3%
}\sup_y||T_{\frac y2}(\int_0^{y}|\varphi_s|^2)sds||^{\frac12}_\infty
||G(f)||_1^{\frac 12}||S(f)||_1^{\frac12}.\qed
\end{eqnarray*}

\medskip \textbf{Proof of Theorem 4.1.} Since
\begin{eqnarray*}
\tau f\varphi^*=4\tau\int_0^\infty \frac {\partial T_{s}f}{\partial
s}\frac {\partial T_{s}\varphi^*}{\partial s}sds =4\tau\int_0^\infty
\frac {\partial T_{3s}f}{\partial s}\frac {\partial
T_{3s}\varphi^*}{\partial s}sds.
\end{eqnarray*}
Setting $\varphi_s=\frac {\partial T_{3s}\varphi}{\partial s}$ and
applying Lemma \ref{fphis}, we get
\begin{eqnarray}
|\tau f\varphi^*|\leq 12\sup_y||T_{\frac y2}(\int_0^{y}|\frac
{\partial T_{3s}\varphi}{\partial s}|^2)sds||^{\frac12}_\infty
||G(f)||_1^{\frac 12}||S(f)||_1^{\frac12}.  \label{3sgfsf}
\end{eqnarray}
On the other hand, we have
\begin{eqnarray}
&&||T_{\frac y2}(\int_0^{y}|\frac {\partial T_{s}\varphi}{\partial
s}|^2)sds||^{\frac12}_\infty  \nonumber \\
&\leq&||T_{\frac y2}(\int_0^{\frac y2}|\frac {\partial
T_{s}\varphi}{\partial s}|^2)sds||^{\frac12}_\infty+ ||T_{\frac
y2}(\int_{\frac y2}^{y}|\frac {\partial T_{s}\varphi}{\partial
s}|^2)sds||^{\frac12}_\infty  \nonumber \\
&\leq&||T_{\frac y2}(\int_0^{\frac y2}|\frac {\partial
T_{s}\varphi}{\partial s}|^2)sds||^{\frac12}_\infty+ ||T_{\frac
y2}(\int_{\frac y2}^{y}T_{\frac y4}|\frac {\partial T_{s-\frac
y4}\varphi}{\partial s}|^2)sds||^{\frac12}_\infty  \nonumber \\
&\leq&||T_{\frac y2}(\int_0^{\frac y2}|\frac {\partial
T_{s}\varphi}{\partial s}|^2)sds||^{\frac12}_\infty+ ||T_{\frac
{3y}4}(\int_{\frac y4}^{\frac {3y}4}|\frac {\partial T_{u}\varphi}{\partial
u}|^2)(u+ \frac y4)du||^{\frac12}_\infty  \nonumber \\
&\leq&||T_{\frac y2}(\int_0^{\frac y2}|\frac {\partial
T_{s}\varphi}{\partial s}|^2)sds||^{\frac12}_\infty+ ||2T_{\frac
{3y}4}(\int_{\frac y4}^{\frac {3y}4}|\frac {\partial T_{u}\varphi}{\partial
u}|^2)udu||^{\frac12}_\infty  \nonumber \\
&\leq&(1+\sqrt2)||\varphi||_{\mathrm{BMO}_c^C}.  \label{trick}
\end{eqnarray}
Using the same idea, we can get
\begin{eqnarray}
||T_{\frac y2}(\int_0^{y}|\frac {\partial T_{3s}\varphi}{\partial
s}|^2)sds||^{\frac12}_\infty =||T_{\frac y2}(\int_0^{3y}|\frac {\partial
T_{s}\varphi}{\partial s}|^2)sds||^{\frac12}_\infty\leq c||\varphi||_{%
\mathrm{BMO}_c^C}.  \label{3s}
\end{eqnarray}
By (\ref{3sgfsf}) and (\ref{3s}), we get
\begin{eqnarray*}
|\tau f\varphi^*|\leq
c||\varphi||_{\mathrm{BMO}_c^C}||G(f)||_1^{\frac
12}||S(f)||_1^{\frac12} \leq c||\varphi||_{\mathrm{BMO}_c^C}||f||_{\mathcal{H%
}_{c,1}^S}. \qed
\end{eqnarray*}

\begin{theorem}
\label{end} Suppose $(T_s)_s$ satisfy the $L^{\frac 12}$ condition (\ref
{Lhalf}) and $T_{2s}\leq cT_s,$ for all $s$ or $T_s\leq cT_{2s}$ for all $s$%
. Then
\[
||f||_{\mathcal{H}_{c,1}^S}\approx ||f||_{\mathcal{H}_{c,1}^G}.
\]
\end{theorem}

\textbf{Proof.} As mentioned in Remark \ref{weakmonotone}, the assumption of
Theorem \ref{end} is sufficient for $(\mathcal{T}_1^{(T_{4s})})^*\subset%
\mathcal{T}^{(T_{4s})}_\infty$.
%Here ${\cal T}_1^{(T_{4s})},{\cal T}^{(T_{4s})}_\infty$ are the tent space associated with $(T_{4s})_s$.
Then
\begin{eqnarray*}
||f||_{\mathcal{H}_{c,1}^S} &=&||(s\frac {\partial
T_{4s}f}{4\partial
s})_s||_{\mathcal{T}_1^{(T_{4s})}} \\
&\leq&c\sup_{||(\varphi_s)_s||_{\mathcal{T}_\infty^{(T_{4s})}\leq1}}{\tau}%
\int_0^\infty \frac {\partial T_{4s}f}{4\partial s}\varphi_ssds \\
&=&c\sup_{||\varphi_s||_{\mathcal{T}^{(T_{4s})}_\infty\leq1}}{\tau}%
\int_0^\infty \frac {\partial T_{3s}f}{3\partial s}T_{s}(\varphi_s)sds \\
&=&\frac {c}3\sup_{||\varphi_s||_{\mathcal{T}^{(T_{4s})}_\infty\leq1}}{\tau}%
\int_0^\infty \frac {\partial T_{3s}f}{\partial s}T_{s}(\varphi_s)sds.
\end{eqnarray*}
By Lemma \ref{fphis}, we get
\begin{eqnarray*}
||f||_{\mathcal{H}_{c,1}^S}&\leq&c\sup_{||\varphi_s||_{\mathcal{T}%
_\infty^{(T_{4s})}\leq1}} \sup_y||T_{\frac
y2}(\int_0^{y}|T_{s}\varphi_s|^2)sds||^{\frac12}_\infty ||G(f)||_1^{\frac
12}||S(f)||_1^{\frac12}.
\end{eqnarray*}
By the assumption $T_{2s}\leq cT_s$ (or $T_s\leq cT_{2s}$) and similar trick
used in (\ref{trick}), we can get
\begin{eqnarray}
\sup_y||T_{\frac y2}(\int_0^{y}|T_{s}\varphi_s|^2)sds||^{\frac12}_\infty\leq
c ||(\varphi_s)_s||_{\mathcal{T}_\infty^{(T_{4s})}}.
\end{eqnarray}
Therefore,
\begin{eqnarray*}
||f||_{\mathcal{H}_{c,1}^S}\leq c||G(f)||_1^{\frac 12}||S(f)||_1^{\frac12}.
\end{eqnarray*}
And
\begin{eqnarray*}
||f||_{\mathcal{H}_{c,1}^S}\leq c||G(f)||_1=c||f||_{\mathcal{H}_{c,1}^G}.
\end{eqnarray*}
The inverse relation is (\ref{HGS}). We then finished the proof.\qed

\section*{Appendix}
We will prove that a large class of semigroups on $\Bbb{R}^n$
(including classical heat semigroup) satisfies the $L^{\frac12}$
condition (\ref{Lhalf}).

\begin{proposition}
\label{thalf} Let $(T_t)_t$ be a semigroup on $\Bbb{R}^n$ with kernel $%
K_t(x,s)$, i.e. $T_t(f)(x)=\int_{\Bbb{R}^n}K_t(x,s)f(s)ds$. Suppose
that there exist constants $r>1,c>0\in \Bbb{R}$ such that
\begin{eqnarray}
K_t(x,s)\leq \frac{c\phi (t)^r}{\phi (t)^{n+r}+|x-s|^{n+r}},  \label{kx}
\end{eqnarray}
with $\phi (t)$ a positive function of $t$. Then $(T_t)_t$ satisfies the $%
L^{\frac 12}$ condition (\ref{Lhalf}).
\end{proposition}

\textbf{Proof.} Fix $t>0$. Let $n=1$. Consider two increasing
filtrations of $\sigma -$algebras: ${\cal D}=\{{D}_k\}_{k\in
\mathbb{Z}}$, with ${D}_k$ the $\sigma -$ algebra generated by the
atoms
\begin{eqnarray*}
D_k^j=({\phi(t)}j4^{-k},{\phi(t)}(j+1)4^{-k}];\quad j\in \mathbb{Z},
\end{eqnarray*}
and ${\cal D}^{\prime}=\{{D}^{\prime}_k\}_{k\in \mathbb{Z}}$, with
${D}^{\prime}_k$ generated by the atoms
\begin{eqnarray*}
D_k^{{\prime }j}=({\phi(t)}(j+\frac 13)4^{-k},{\phi(t)}(j+\frac
43)4^{-k}],\quad j\in \mathbb{Z}.
\end{eqnarray*}
%It is easy to verify that ${\mathcal D}^{\prime }=\{{\mathcal
%D}_k^{\prime }\}_{n\in \mathbb{Z}}$ is indeed an increasing
%filtration.
Let $E_k,E_k^{\prime}$ be the conditional expectation with respect to $D_k$
and $D_k^{\prime}$. It is easy to verify that, for any $f\geq0$,
%\begin{proposition} (see [Mei], Proposition 3.1)
%For any interval $I\subset \Bbb{R},$ there exist $k_I,N\in \Bbb{Z}$
%such
%that $I\subset $ $D_N^{k_I}$ and $|D_N^{k_I}|\leq 6|I|$ or $I\subset $ $%
%D_N^{^{\prime }k_I}$ and $|D_N^{^{\prime }k_I}|\leq 6|I|,$ the
%constant $N$ only depends on the length of $I.$
%\end{proposition}
\begin{eqnarray}
E_k(f)&\leq& 4E_{k-1}(f),\ \ E_k^{\prime}\leq 4E_{k-1}^{\prime}(f); \\
E_k(f)&\leq& 3E_k^{\prime}E_k(f),\ \ E_k^{\prime}\leq 3E_kE_k^{\prime}(f); \\
T_t(f)&\leq&
c\sum_{k=-\infty}^04^{kr}E_k(f)+c\sum_{k=-\infty}^04^{kr}E^{\prime}_k(f)
\end{eqnarray}
Therefore, for any $f,g\geq0$,
\begin{eqnarray*}
&&T_t(fT_tg) \\
&\leq&
c\sum_{k,i=-\infty}^04^{kr}4^{ir}[E_k(fE_ig)+E^{\prime}_k(fE_i^{%
\prime}g)+E_k^{\prime}(fE_ig)+E_k(fE^{\prime}_ig)] \\
&\leq& c\sum_{k\geq
i}4^{kr}4^{ir}(E_kfE_ig+E^{\prime}_kfE^{\prime}_ig)+c%
\sum_{k<i}4^{kr}4^{ir}4^{i-k}[E^{\prime}_k(fE_kg)+E_k(fE^{\prime}_kg)] \\
&\leq& c\sum_{k\geq
i}4^{kr}4^{ir}(E_kfE_ig+E^{\prime}_kfE^{\prime}_ig)+3c%
\sum_{k<i}4^{kr}4^{ir}4^{i-k}[E^{\prime}_k(fE_k^{\prime}E_kg)+E_k(fE_kE^{%
\prime}_kg)] \\
&\leq& c\sum_{k\geq
i}4^{kr}4^{ir}(E_kfE_ig+E^{\prime}_kfE^{\prime}_ig)+3c%
\sum_{k<i}4^{kr}4^{ir}4^{i-k}[(E^{\prime}_kf)(E_k^{%
\prime}E_kg)+(E_kf)(E_kE^{\prime}_kg)].
\end{eqnarray*}
And
\begin{eqnarray*}
\int_{\Bbb{R}}[T_t(fT_tg)]^{\frac12}&\leq& c\sum_{k\geq i}2^{kr}2^{ir}\int_{%
\Bbb{R}}[(E_kf)^\frac12(E_ig)^\frac12+(E^{\prime}_kf)^\frac12(E^{%
\prime}_ig)^\frac12] \\
&&+c\sum_{k<i}2^{kr}2^{ir}2^{i-k}\int_{\Bbb{R}}[(E^{\prime}_kf)^{\frac12}
(E_k^{\prime}E_kg)^\frac12+(E_kf)^\frac12(E_kE^{\prime}_kg)^\frac12] \\
&\leq& c\sum_{k\geq
i}2^{kr}2^{ir}2||f||_1^\frac12||g||_1^\frac12+c%
\sum_{k<i}2^{kr}2^{ir}2^{i-k}2||f||_1^\frac12||g||_1^\frac12 \\
&\leq&
c||f||_1^\frac12||g||_1^\frac12+2c\sum_{k<i}2^{k(r-1)}2^{i(r+1)}||f||_1^%
\frac12||g||_1^\frac12 \\
&\leq&
c||f||_1^\frac12||g||_1^\frac12+2c_r\sum_{i=-%
\infty}^02^{i(r-1)}2^{i(r+1)}||f||_1^\frac12||g||_1^\frac12 \\
&\leq& c_r||f||_1^\frac12||g||_1^\frac12.
\end{eqnarray*}
Then $(T_t)_t$ satisfies the $L^\frac12$ condition (\ref{Lhalf}).

For $n>1$, we use the filtrations in Remark 7 of [M] and can prove the
proposition by the same idea presented above.\qed

Since classical heat semigroup on $\Bbb{R}^n$ is a convolution
operator with a kernel
\begin{eqnarray*}
K_t(x)&=&\frac{\exp (-\frac{|x|^2}{4t})}{(4\pi t)^{\frac n2}}
\end{eqnarray*}
which satisfies (\ref{kx}) with $\phi(t)=2t^\frac12$. We then get

\begin{corollary}\label{lastcor}
Classical heat semigroup $(T_t)_t$ satisfies the $L^{\frac 12}$
condition (\ref {Lhalf}).
\end{corollary}

\begin{remark}{\rm
Another way to prove Corollary \ref{lastcor} is to verify the
condition (i) of Remark \ref{equivalent}. The proof will be indirect
but easier and will imply that $(T_t\otimes I)_t$ satisfies the
$L^{\frac 12}$ condition as well on $L^{\infty}({\Bbb R}^n)\otimes
B(\ell_2)$ with $I$ the identity operator on $B(\ell_2)$.

In a forth coming paper with Avsec Stephen, we are going to use this
property of $(T_t\otimes I)_t$ to prove an $H^1$-BMO duality result
on group von Neumann algebra $VN(G)$. The idea is to embed $VN(G)$
into the crossed product $L^\infty({\Bbb R}^n)\rtimes G$.
}\end{remark}

\medskip
{\bf Acknowledgment.} The author is grateful to M. Junge and Q. Xu
for helpful discussions. The author also thanks the referee for a
careful reading and useful comments. The author thanks the
organizers of the workshop in Analysis and Probability in College
Station, Tx, where part of this work was carried out.

\textbf{Reference}

%[AC] C. Anantharaman-Delaroche, On ergodic theorems for free group actions
%on noncommutative spaces. Probab. Theory Related Fields 135 (2006), no. 4,
%520--546.

%[BL] J. Bergh, J. L\"{o}fstr\"{o}m, Interpolation spaces. An introduction.
%Grundlehren der Mathematischen Wissenschaften, No. 223. Springer-Verlag,
%Berlin-New York, 1976.

[CMS] Coifman, R. R.; Meyer, Y.; Stein, E. M. Some new function spaces and
their applications to harmonic analysis. J. Funct. Anal. 62 (1985), no. 2,
304--335.

[Da] E.B. Davies, Heat kernels and spectral theory, Cambridge Univ. Press,
1989.

[DY] X. Duong, L. Yan, Duality of Hardy and BMO spaces associated with
operators with heat kernel bounds. J. Amer. Math. Soc. 18 (2005), no. 4,
943--973.

[G]  J. B. Garnett, Bounded Analytic Functions, Pure and Applied
Mathematics, 96. Academic press, Inc., New York-London, 1981.

%[HM] S. Hofmann, S Mayboroda, Hardy and BMO spaces associated to divergence form elliptic operators, Math. Ann., to appear.

[J1] M. Junge, Doob's Inequality for Non-commutative Martingales, J. Reine
Angew. Math. 549 (2002), 149-190.

[J2] M. Junge, Square function and Riesz-transform estimates for
subordinated semigroups, preprint.

[JLX] M. Junge, C. Le Merdy, Q. Xu, $H^\infty $ Functional calculus
and square functions on noncommutative $L^p$ -spaces. Ast\'{e}risque
No. 305 (2006), vi+138 pp.

[JM] M. Junge, T. Mei, Noncommutative Riesz Transforms-A Probabilistic
Approach, Preprint.

[JX] M. Junge, Q. Xu, Noncommutative maximal ergodic theorems. J. Amer.
Math. Soc. 20 (2007), no. 2, 385--439.

[JS] M. Junge, D. Sherman, Noncommutative $L^p$ modules, J. of
operator theory, .

%[JZ] M. Junge, Z. J. Ruan, Decomposable maps on non-commutative $L\sb p$%
%-spaces. Operator algebras, quantization, and noncommutative geometry,
%355--381, Contemp. Math., 365, Amer. Math. Soc., Providence, RI, 2004.

[La] E. C. Lance, Hilbert $C^*$- Modules: A Toolkit for Operator
Algebraists, London Math. Soc. Lecture Note Ser., vol. 210,
Cambridge Univ. Press, Cambridge 1995.

%[Lp1] F. Lust-Piquard, Riesz transforms associated with the number operator
%on the Walsh system and the fermions. J. Funct. Anal. 155 (1998), no. 1,
%263--285.

%[Lp2] F. Lust-Piquard, Riesz transforms on deformed Fock spaces. Comm. Math.
%Phys. 205 (1999), no. 3, 519--549.

%[Lp3] F. Lust-Piquard, Dimension free estimates for discrete Riesz
%transforms on products of abelian groups. Adv. Math. 185 (2004), no. 2,
%289--327.

[MTX] T. Mart\'{i}nez, J. L. Torrea, Q. Xu, Vector-valued
Littlewood-Paley-Stein theory for semigroups. Adv. Math. 203 (2006),
no. 2, 430--475.

[M] T. Mei, Operator Valued Hardy Spaces, Memoirs of AMS, 2007, V.
188, No. 881.

[Ne] E. Nelson, Notes on non-commutative integration, J. Funct. Anal., 15
(1974), 103-116.

[Pa] W. Paschke, Inner product modules over $B^*$ algebras,  Tran.
of AMS, 182 (1973) 443-468.

[Pe] K. Petersen, Brownian motion, Hardy spaces and Bounded Mean
Oscillation, LMS Lecture Notes Series 28, Cambridge University
Press, 1977.

[P1] G. Pisier, Non-commutative vector valued $L_ p$-spaces and completely $%
p $-summing maps. Ast\'{e}risque No. 247 (1998), vi+131 pp.

%[Mey] P. A. Meyer, D\'{e}monstration probabiliste de certaines
%in\'{e}galit\'{e}s de Littlewood-Paley. IV. Semi-groupes de convolution
%sym\'{e}triques. (French) S\'{e}minaire de Probabilit\'{e}s, X (Premi\`{e}re
%partie, Univ. Strasbourg, Strasbourg, ann\'{e}e universitaire 1974/1975),
%pp. 175--183. Lecture Notes in Math., Vol. 511, Springer, Berlin, 1976.

[PX] G. Pisier, Q. Xu, Non-commutative $L_p$-spaces. Handbook of the
geometry of Banach spaces, Vol. 2, 1459--1517, North-Holland, Amsterdam,
2003.

[PX1] G. Pisier and Q. Xu, Non-commutative martingale inequalities.
Comm. Math. Phys. Vol. 189, (1997), 667-698.

[St1] E. M. Stein, Harmonic Analysis, Princeton Univ. Press, Princeton, New
Jersey, 1993.

[St2] E. M. Stein, Topic in Harmonic Analysis (related to Littlewood-Paley
theory), Princeton Univ. Press, Princeton, New Jersey, 1970.

[V1] N. Varopoulos, Aspects of probabilistic Littlewood-Paley
theory. J. Funct. Anal. 38 (1980), no. 1, 25--60.

%\bigskip $
%\begin{array}{l}
%\mbox{Math. Dept.} \\
%\mbox{Univ. of Illinois at Urbana Champaign} \\
%\mbox{Urbana, IL, 61801} \\
%\mbox{U. S. A.} \\
%\mbox{mei@math.uiuc.edu}
%\end{array}
%$

\end{document}